\documentclass[11pt,twoside]{article} 
\usepackage{latexsym} 
\input{amssym.def} 
\input{amssym} 
\newfont{\fra}{eufm10 scaled 1095} 
\newfont{\Bb}{msbm10 scaled 1095} 
\newfont{\Bbg}{msbm10 scaled 1680} 
\newcommand\CC{{\mbox{\Bb C}}} 
\newcommand\RR{{\mbox{\Bb R}}} 
\newcommand\NN{{\mbox{\Bb N}}} 
\newcommand\ZZ{{\mbox{\Bb Z}}}

\newcommand\fb{{\frak b}} 
\newcommand\fg{{\frak{g}}} 
\newcommand\fh{{\frak h}} 
\newcommand\fri{{\frak i}} 
 
\newcommand\fl{{\frak l}} 
 
\newcommand\fn{{\frak n}} 
\newcommand\fk{{\frak k}}

\newcommand\fa{{\frak a}} 
\newcommand\fd{{\frak d}} 
\newcommand\fr{{\frak r}} 
\newcommand\fs{{\frak s}} 
\newcommand\fv{{\frak v}} 
\newcommand\fz{{\frak z}} 

\newcommand\cH{{{\cal H}}} 
\newcommand\cA{{\cal A}} 
\newcommand\cL{{\cal L}}

\newcommand\ph{\varphi}

\newcommand{\fsl}{\mathop{{\frak s \frak l}}}

\newcommand{\fso}{\mathop{{\frak s \frak o}}} 
 
\newcommand{\fsu}{\mathop{{\frak s \frak u}}} 

\newcommand{\Der}{\mathop{{\rm Der}}} 
\newcommand{\Iso}{\mathop{{\rm Iso}}} 
\newcommand{\GL}{\mathop{{\rm GL}}} 
\newcommand{\grO}{{\rm O}}
\newcommand{\Hom}{\mathop{{\rm Hom}}} 
\newcommand{\Aut}{\mathop{{\rm Aut}}}  
\newcommand{\Out}{\mathop{{\rm Out}}}
\newcommand{\Sym}{{{\rm Sym}}}
\newcommand{\Id}{{{\rm id}}} 
\newcommand{\ad}{{{\rm ad}}}

\newcommand{\Ker}{\mathop{{\rm ker}}}

\renewcommand{\Re}{\mathop{{\rm Re}}} 
\renewcommand{\Im}{\mathop{{\rm Im}}} 
\newcommand{\diag}{\mathop{{\rm diag}}}

\newcommand{\Span}{{{\rm span}}} 
\newcommand{\mod}{\mathop{{\rm mod}}} 
\newcommand{\proj}{{{\rm pr}}} 
\newcommand\ip{{\langle\cdot \,,\cdot \rangle}} 
\newcommand\lb{{[\cdot\,,\cdot]}} 
 
\newcommand\dd{\fd_{\alpha,\gamma}(\fl,\Phi_\fl,\fa)} 
\newcommand\HPhi{{\cH^2_Q(\fl,\Phi_\fl,\fa)}} 
\newcommand\ZPhi{{{\cal Z}^2_Q(\fl,\Phi_\fl,\fa)}}  
\newcommand\proof{{\sl Proof. }} 
\newcommand{\qed}{\hspace*{\fill}\hbox{$\Box$}\vspace{2ex}} 
\newtheorem{theo}{Theorem}[section] 
\newtheorem{pr}[theo]{Proposition}
\newtheorem{de}[theo]{Definition}
\newtheorem{ex}[theo]{Example}
\newtheorem{re}[theo]{Remark}
\newtheorem{co}[theo]{Corollary}
\newtheorem{lm}[theo]{Lemma}
\begin{document} 
\title{Indefinite extrinsic symmetric spaces II} 
\author{Ines Kath}
\maketitle 
\begin{abstract}
\noindent Having developed a description of indefinite extrinsic symmetric spaces by corresponding infinitesimal objects in the preceding paper  we now study the classification problem for these algebraic objects.  In most cases the transvection group of an indefinite extrinsic symmetric space is not semisimple, which makes the classification difficult. We use the recently developed method of quadratic extensions for $(\fh,K)$-invariant metric Lie algebras to tackle this problem. We obtain a  one-to-one correspondence between isometry classes of extrinsic symmetric spaces and a certain cohomology set.  This allows a systematic construction of extrinsic symmetric spaces and explicit classification results, e.g., if the metric of the embedded manifold or the ambient space has a small index. We will illustrate this by classifying all Lorentzian extrinsic symmetric spaces.\\[2ex]

\end{abstract}
\section{Introduction} 
We study special embeddings of symmetric spaces in pseudo-Euclidean spaces. 
Recall that a non-degenerate submanifold $M$ of a pseudo-Euclidean space $V$ is called extrinsic symmetric if it is invariant under the reflection at each of its affine normal spaces. In particular, $M$ is an abstract symmetric space. Extrinsic symmetric spaces can also be characterised as those connected complete 
submanifolds whose second fundamental form is parallel. 
Two extrinsic symmetric spaces $M\hookrightarrow V$ and $M'\hookrightarrow V'$ are called isometric, if there is an isometry $V\rightarrow V'$ that maps $M$ to $M'$. In the present paper we continue to discuss the classification of extrinsic symmetric spaces up to isometry. The first step was done in \cite{Kext1}, were we gave an algebraic description of extrinsic symmetric spaces.  Now we want to classify these algebraic objects.

Let us shortly summarise the algebraic description developed in \cite{Kext1}. It generalises Ferus' construction \cite{F2,F3,EH} to the pseudo-Riemannian situation.  An extrinsic symmetric space $M\hookrightarrow V$ is called full if it is not contained in a proper subspace of $V$. We proved that full extrinsic symmetric spaces correspond to so-called full  extrinsic symmetric triples. Roughly speaking, an extrinsic symmetric triple $(\fg,\Phi,\ip)$ consists of a Lie algebra $\fg$, a (non-degenerate) $\fg$-invariant scalar product $\ip$ and a pair  $\Phi=(D,\theta)$, where $\theta$ is an isometric involution on $\fg$ and $D$ is an antisymmetric derivation satisfying $D^3=-D$  and certain further conditions, see Section~\ref{S2} for an exact definition.

Obviously, we cannot assume without loss of generality that $M\hookrightarrow V$ is full since the minimal subspace that contains $M$ can be degenerate. Therefore, as an intermediate step, we allow $V$ to be degenerate, while $M$ remains non-degenerate. Then we can confine to full submanifolds. Now full extrinsic symmetric spaces correspond to so-called full weak extrinsic symmetric triples, where weak means that $V_+\subset \fg$ can be degenerate. In \cite{Kext1} we have seen that if $(\tilde\fg,\Phi,\ip)$ is a full weak extrinsic symmetric triple and $R=\tilde\fg\cap\tilde\fg^\perp$ is its metric radical, then $R$ is central and $\fg:=\tilde\fg/R$ is a full extrinsic symmetric triple.  Isomorphism classes of full weak extrinsic symmetric triples that arise in this way as a central extension of $(\fg,\Phi,\ip)$ by a vector space $R$ can be described by a suitable subset of $H^2(\fg,R)/(\Aut(\fg,\Phi,\ip)\times \GL(R))$, see Section \ref{Sfwt} for more details. This reduces the classification of extrinsic symmetric spaces to the following

{\bf Problem:} Determine all full extrinsic symmetric triples  $(\fg,\Phi,\ip)$ up to isomorphisms.

 As in the case of ordinary symmetric spaces the semisimple situation is well studied. Naitoh \cite{Nai} classified all extrinsic symmetric triples $(\fg,\Phi,\ip)$ with semisimple $\fg$. His results are based on  the classification of non-degenerate Jordan triple systems due to Neher \cite{Ne1, Ne2, Ne3}. 
The non-semisimple case is much more involved. Recall that it is impossible to give an explicit classification of all non-semisimple pseudo-Riemannian symmetric spaces. We will see that the same is true for extrinsic symmetric spaces.  Therefore we are looking for methods that allow a systematic study also in this situation. This will be the main goal of the present paper. In the first part of the paper we present a structure theory for non-semisimple extrinsic symmetric triples, which relies on the more general theory for $(\fh,K)$-equivariant metric Lie algebras developed in \cite{KOesi}. 
This will lead to a description of the set of isomorphism classes of extrinsic symmetric triples by means of a certain cohomology set. More exactly, using the method of quadratic extensions we will prove that isomorphism classes of full extrinsic symmetric triples without simple ideals correspond to elements of 
$$
 \coprod  \HPhi_{\sharp}/ G_{\fl,\Phi_{\fl},\fa}\ ,
$$
where the union is taken over a set of representatives of isomorphism
classes of pairs $((\fl,\Phi_\fl),\fa)$ consisting of an $(\RR,\ZZ_2)$-equivariant Lie algebra  $(\fl,\Phi_\fl)$ and a semi-simple orthogonal $(\fl,\Phi_\fl)$-module $\fa$. 
Under suitable additional conditions this cohomology set can be computed, which leads to explicit classification results. A typical example of such a suitable condition is the restriction to small dimensions or to a small index  of the metric of the embedded or the ambient space. 
Furthermore, this description of isomorphism classes gives a method for a systematic construction of examples of extrinsic symmetric spaces.

In the second part of the present paper we will illustrate this method by classifying extrinsic symmetric spaces for which the embedded space has Lorentzian signature. The result is known for surfaces in the four-dimensional Minkowski space \cite{CV}. Here we will consider the general case. We will be mainly interested in indecomposable extrinsic symmetric spaces. Here decomposable means that not only the embedded manifold $M$ is a product of symmetric spaces but also the embedding $M\hookrightarrow V$ itself decomposes into a product of embeddings. Recall that each indecomposable (ordinary) Lorentzian symmetric space is either a space of constant curvature or it has a solvable transvection group and is covered by a Cahen-Wallach space \cite{CW}. The extrinsic symmetric spaces of nonzero constant curvature can be read from Naitoh's list. Besides the natural embeddings of the Lorentzian spheres $S^{1,n-1}$ and the Lorentzian hyperboloids $H^{1,n-1}$ embeddings of products $S^k\times S^{1,l}$ and $H^k\times H^{1,l}$ appear, which are not decomposable as extrinsic symmetric spaces. The classification of flat Lorentzian extrinsic spaces and those which are covered by a Cahen-Wallach space is the difficult part  of our task. To solve it we will apply the structure theory developed in the first part of the paper. Let us first consider  full and indecomposable flat Lorentzian extrinsic symmetric spaces. We will see that minimal and non-minimal ones exist, where as usual a submanifold is called minimal if its mean curvature vector vanishes. As minimal embeddings the identity map $\RR^{1,0}\rightarrow \RR^{1,0}$ and certain embeddings $\RR^{1,1}\hookrightarrow \RR^{1,2}$, $\RR^{1,1}\hookrightarrow \RR^{2,1}$, and $\RR^{1,2}\hookrightarrow \RR^{2,3}$ occur. In the non-minimal case we have embeddings
of flat Lorentzian manifolds of dimension $2+l$ into $\RR^{2,2+2l}$ or into  $\RR^{2+l,2+l}$ for $l\in\NN$.
Now let us turn to the case where $M$ is (covered by) a Cahen-Wallach space. We will see that only (quotients of) Cahen-Wallach spaces with very special parameters can be embedded extrinsically symmetric. This is in contrast to the Riemannian case, where almost all symmetric spaces have  realisations as extrinsic symmetric spaces. The precise classification result for the full case is formulated in Theorem~\ref{Tfull}. For all isomorphism classes of full and indecomposable Lorentzian extrinsic symmetric spaces we will also give explicit realisations. Moreover we will describe all extensions to weak extrinsic symmetric triples, which gives  a classification also in the non-full case. For the result see Theorem~\ref{Tende}.

I would like to thank Jost-Hinrich Eschenburg and Martin Olbrich 
for many helpful discussions.    
         
\section{Extrinsic symmetric triples} \label{S2}

We want to recall the definition of an extrinsic symmetric triple from \cite{Kext1}. Here we will use the language of $(\fh,K)$-equivariant metric Lie algebra in the sense of \cite{KOesi}. Later on this will allow us to apply the classification machinery developed in \cite{KOesi}.  
We do not want to explain the general notion of an $(\fh,K)$-equivariant Lie algebra but confine the definition to the special case $\fh=\RR$, $K=\ZZ_2$. 
Let $\ZZ_2=\{1,-1\}$ act on $\RR$ by multiplication. 
\begin{de} An $(\RR,\ZZ_2)$-equivariant Lie algebra in the sense of {\rm 
\cite{KOesi}} is a pair 
$(\fl,\Phi_\fl)$, where $\fl$ is a Lie algebra and $\Phi_\fl$ is a pair 
$(D_\fl,\theta_\fl)$ 
consisting of a semisimple derivation 
$D_\fl:\fl\rightarrow\fl$ and 
an involution $\theta_\fl:\fl\rightarrow \fl$ 
satisfying 
$D_\fl\circ\theta_\fl=-\theta_\fl\circ D_\fl$. 
\end{de} 
We will say that an $(\RR,\ZZ_2)$-equivariant Lie algebra $(\fl,\Phi)$, $\Phi=(D,\theta)$ is $h$-graded if $D^3=-D$ holds. In this case we define $\tau:=\exp(\pi D)$ and
$$\fl_+:=\{X\in\fl\mid \theta(X)=X\},\ \fl_-:=\{X\in\fl\mid \theta(X)=-X\}$$
and
$$\fl^+:= \{X\in\fl\mid  \tau(X)=X \},\ \fl^-:= \{ X\in\fl\mid \tau(X)=-X \}.$$
Furthermore, 
$$ \fl_+^-:=\fl_+\cap\fl^-,\  \fl_+^+:=\fl_+\cap\fl^+,\ \fl_-^-:= \fl_-\cap\fl^-,\  \fl^+_-:= \fl_-\cap\fl^+.$$
Then it is easy to show that 
\begin{equation}\label{ED}
D|_{\fl^+}=0,\ D(\fl_+^-)= \fl_-^-,\ (D|_{\fl^-})^2=-\Id.
\end{equation}
\begin{de} 
An $(\RR,\ZZ_2)$-equivariant metric Lie algebra is a triple $(\fg,\Phi,\ip)$ 
such that 
\begin{itemize} 
\item[(i)] $(\fg,\ip)$ is a metric Lie algebra, 
\item[(ii)] $(\fg,\Phi)$ is an $(\RR,\ZZ_2)$-equivariant Lie algebra, 
\item[(iii)] $D$ is antisymmetric and $\theta$ is an isometry. 
\end{itemize} 
\end{de} 

In particular, $D|_{\fg_+^-}:\fg_+^-\rightarrow \fg_-^-$ is an isometry.

\begin{de} 
An extrinsic symmetric triple is an $(\RR,\ZZ_2)$-equivariant metric Lie algebra $(\fg,\Phi,\ip)$ that is $h$-graded and satisfies  $[\fg_+^-,\fg_+^-]=\fg_+^+$. 
\end{de}
An extrinsic symmetric triple is called full if it satisfies $[\fg_+^-,\fg_-^-]=\fg_-^+$. We can formulate this equivalently in the following way.

\begin{lm}\label{Lfest}
A full extrinsic symmetric triple is an $h$-graded $(\RR,\ZZ_2)$-equivariant metric Lie algebra $(\fg,\Phi,\ip)$ that satisfies $[\fg^-,\fg^-]=\fg^+$. 
\end{lm}
\proof We have to prove that the conditions $[\fg^-_+,\fg^-_+]=\fg^+_+$ and $[\fg^-_+,\fg^-_-]=\fg^+_-$ together are equivalent to $[\fg^-,\fg^-]=\fg^+$. Obviously, $[\fg^-_+,\fg^-_+]=\fg^+_+$ and $[\fg^-_+,\fg^-_-]=\fg^+_-$ imply $[\fg^-,\fg^-]=\fg^+$. Now assume that $[\fg^-,\fg^-]=\fg^+$ holds. Note that (\ref{ED}) implies $[\fg_+^-,\fg_+^-]=[D\fg_-^-,D\fg_-^-]=[\fg_-^-,\fg_-^-]$. Hence it follows from $[\fg^-,\fg^-]=\fg^+$ that  $\fg_+^+=[\fg_+^-,\fg_+^-]+[\fg_-^-,\fg_-^-]=[\fg_+^-,\fg_+^-]$. Moreover, $[\fg^-,\fg^-]=\fg^+$ yields $[\fg^-_+,\fg^-_-]=\fg^+_-$, hence the triple is full. 
\qed

Weak extrinsic symmetric triples $(\fg,\Phi,\ip)$ are defined in the same way as extrinsic symmetric triples, except that $\ip$ is allowed to be degenerate on $\fg_-^+$. 
 
A (weak) extrinsic symmetric triple  is called decomposable if it is the direct sum of two non-trivial (weak) extrinsic symmetric triples. Otherwise it is called indecomposable. An extrinsic symmetric triple decomposes if and only if the associated extrinsic symmetric space $M\hookrightarrow V$ is decomposable, i.e., if it decomposes into a product of non-trivial embeddings $M_1\hookrightarrow V_1$ and $M_2\hookrightarrow V_2$, where $V=V_1\oplus V_2$. Obviously, to reach our aim of a classification of all full extrinsic symmetric triples it will be sufficient to consider  indecomposable ones. 

\begin{de} 
Let $(\fl,\Phi_\fl)$ be an $(\RR,\ZZ_2)$-equivariant Lie algebra. An orthogonal 
$(\fl,\Phi_\fl)$-module 
$(\rho,\fa,\ip_\fa,\Phi_\fa)$ (often abbreviated to $\fa$) consists of 
\begin{itemize} 
\item a pseudo-Euclidean vector space $(\fa,\ip_\fa)$, 
\item a pair $\Phi_\fa=(D_\fa,\theta_\fa)$, where $\theta_\fa\in 
O(\fa,\ip_\fa)$ is an involution and 
$D_\fa\in\fso(\fa,\ip_\fa)$ is a semisimple map satisfying $D_\fa\circ\theta_\fa=-\theta_\fa\circ  D_\fa$, 
\item an orthogonal representation $\rho:\fl\rightarrow \fso(\fa,\ip)$ that 
satisfies 
\begin{equation}\label{EDrho}\rho(\theta_\fl(L))=\theta_\fa\circ\rho(L)\circ\theta_\fa^{-1},\ 
\rho(D_\fl(L))=[D_\fa,\rho(L)]
\end{equation} 
for all $L\in\fl$. 
\end{itemize} 
We will say that $\fa$ is $h$-graded if $D_\fa^3=-D_\fa$. In this case we define $\fa^+,\fa^-,\fa_+,\fa_-$ etc.~as usual.
\end{de} 
Homomorphisms of (metric) $(\RR,\ZZ_2)$-equivariant Lie 
algebras are defined in the obvious way. 
Let $(\fl_{i},\Phi_{\fl_{i}})$, $i=1,2$, be two $(\RR,\ZZ_2)$-equivariant Lie 
algebras and let $(\rho_{i},\fa_{i})$, $i=1,2$, be orthogonal 
$(\fl_{i},\Phi_{\fl_{i}})$-modules. Let $S:\fl_1\rightarrow \fl_{2}$ 
be a homomorphism of $(\RR,\ZZ_2)$-equivariant Lie algebras and let 
$U:\fa_{2}\rightarrow \fa_{1}$ be an $(\RR,\ZZ_2)$-equivariant isometric embedding. 
Suppose that 
$$U\circ\rho_{2}(S(L))=\rho_{1}(L)\circ U$$ holds for all $L\in\fl$. 
Then we call $(S,U)$ a morphism of pairs. We will write this as 
$(S,U):((\fl_{1},\Phi_{\fl_{1}}),\fa_{1})   \rightarrow 
((\fl_{2},\Phi_{\fl_{2}}),\fa_{2})$.

\section{Quadratic cohomology and quadratic extensions}
\label{S3}
In this section we recall some facts on quadratic cohomology from \cite{KOesi}, where we specialise again to the case $(\fh,K)=(\RR,\ZZ_2)$.
Let $(\fl,\Phi_\fl)$ be an $(\RR,\ZZ_2)$-equivariant Lie algebra and $\fa$ an orthogonal $(\fl,\Phi_\fl)$-module. As usual, $(C^*(\fl,\fa),d)$ denotes the standard complex for Lie algebra cohomology with coefficients in $\fa$. If $\fa=\RR$ is the trivial $(\fl,\Phi_\fl)$-module, then we abbreviate this to $(C^*(\fl),d)$.
Let $C^p(\fl,\fa)^{(D,\theta)}$ be the space of cocycles in $C^p(\fl,\fa)$ that are invariant under $D$ and $\theta$. In the general case of an $(\fh,K)$-action on $\fa$ and $\fl$ considered in \cite{KOesi} this space is denoted by $C^p(\fl,\fa)^{(\fh,K)}$. The set of quadratic $1-$cochains is defined as
$${\cal C} ^1_Q(\fl,\Phi_{\fl},\fa)= C^1(\fl,\fa)^{(D,\theta)}\oplus 
C^2(\fl)^{(D,\theta)}. 
$$ 
This set is a group with group operation defined by 
$$(\tau_1,\sigma_1)*(\tau_2,\sigma_2)=(\tau_1+\tau_2, \sigma_1 
+\sigma_2 +\textstyle{\frac12} \langle \tau_1\wedge 
\tau_2\rangle)\,,$$ 
where $\langle\cdot\,,\cdot\rangle$ denotes the usual wedge product followed by the inner product on $\fa$.
Let us consider the set
$${\cal Z} ^{2}_Q(\fl,\Phi_{\fl},\fa)=\{(\alpha,\gamma) \in 
C^{2}(\fl,\fa)^{(D,\theta)}\oplus 
C^{3}(\fl)^{(D,\theta)} \mid d\alpha=0,\ 
d\gamma=\textstyle{\frac12}\langle\alpha \wedge\alpha\rangle\} $$ 
whose elements are called quadratic $2$-cocycles. The group
${\cal C} ^{1}_Q(\fl,\Phi_{\fl},\fa)$ acts from the right on ${\cal Z} 
^{2}_Q(\fl,\Phi_{\fl},\fa)$ by 
$$(\alpha,\gamma)(\tau,\sigma)=\left(\,\alpha +d\tau,\gamma +d\sigma 
+\langle(\alpha +\textstyle{\frac12} d\tau)\wedge\tau\rangle\,\right).$$ 
The quadratic cohomology set is defined as the orbit space 
$${\cal H}^{2}_Q(\fl,\Phi_{\fl},\fa):={\cal Z}^{2}_Q(\fl,\Phi_{\fl},\fa)/ {\cal 
C}^{1}_Q(\fl,\Phi_{\fl},\fa).$$ 
The equivalence class of 
$(\alpha,\gamma)\in {\cal Z} ^{2}_Q(\fl,\Phi_{\fl},\fa)$ in ${\cal 
H}^{2}_Q(\fl,\Phi_{\fl},\fa)$ is denoted by $[\alpha,\gamma]$. 

Each morphism of pairs $(S,U):((\fl_{1},\Phi_{\fl_{1}}),\fa_{1})   \rightarrow ((\fl_{2},\Phi_{\fl_{2}}),\fa_{2})$ induces maps
\begin{eqnarray*} 
&&(S,U)^*:C^k(\fl_2,\fa_2) \rightarrow C^k(\fl_1,\fa_1),\  (S,U)^*(\omega)= U\circ S^*\omega\\
&&(S,U)^*: {\cal H}^{2}_Q(\fl_2,\Phi_{\fl_2},\fa_2)\rightarrow {\cal 
H}^{2}_Q(\fl_1,\Phi_{\fl_1},\fa_1), \ (S,U)^*(\alpha,\gamma):=(U\circ S^*\alpha, S^*\gamma).
\end{eqnarray*}

Given an $(\RR,\ZZ_2)$-equivariant Lie algebra $(\fl,\Phi_\fl)$, an orthogonal $(\fl,\Phi_\fl)$-module $\fa$ and a quadratic 2-cocycle $(\alpha,\gamma)$ we now construct an $(\RR,\ZZ_2)$-equivariant metric Lie algebra.
We consider the vector space 
$\fd:=\fl^*\oplus\fa\oplus \fl$ 
and define an inner product $\ip$ on $\fd$ by 
$$ 
\langle Z_1+A_1+L_1, Z_2+A_2+L_2\rangle :=   
\langle A_1,A_2\rangle _\fa +Z_1(L_2)+Z_2(L_1) 
$$ 
for $Z_1, Z_2\in \fl^*,\, A_1,A_2\in\fa$ and $L_1,L_2\in\fl$. The 
antisymmetric bilinear map $\lb:\fd\times \fd\rightarrow\fd$ defined by 
$[\fl^*,\fl^*\oplus\fa]=0$ and 
\begin{eqnarray*} 
\,[L_1,L_2]&=&\gamma(L_1,L_2,\cdot)+\alpha(L_1,L_2)+[L_1,L_2]_\fl\\ 
\,[L,A] &=& -\langle A,\alpha(L,\cdot)\rangle +L(A)\\ 
\,[A_1,A_2] &=& \langle\rho(\cdot) A_1,A_2\rangle \\ 
\,[L,Z] &=& \ad^*(L)(Z) 
\end{eqnarray*} 
for $L,L_1,L_2\in\fl$, $A,A_1,A_2\in\fa$ and $Z\in\fl^*$ is a Lie bracket. The scalar product $\ip$ is invariant 
with respect to this bracket \cite{KO1}.
Moreover, we define a pair $\Phi=(D,\theta)$  by
\begin{eqnarray*} 
D(Z+A+L) &=& -D_\fl^{*}(Z)+D_{\fa}(A)+D_{\fl}(L)\\ 
\theta(Z+A+L) &=&\theta_{\fl}^{*}(Z)+\theta_{\fa}(A)+\theta_{\fl}(L). 
\end{eqnarray*} 
Then $\fd_{\alpha,\gamma}(\fl,\Phi_{\fl},\fa):=(\fd,\Phi,\ip)$ is an 
$(\RR,\ZZ_2)$-equivariant metric Lie algebra. It is called quadratic extension of $(\fl,\Phi_\fl)$ by $\fa$ associated with the cocycle $(\alpha,\gamma)\in\ZPhi$.
Obviously, $\dd$ is $h$-graded if and only if $(\fl,\Phi_\fl)$ and $\fa$ are so.

\section{Structure theory for extrinsic symmetric triples} \label{Sst}
First we will show that we can decompose any extrinsic symmetric triple into a semisimple one and one that does not have simple ideals. As remarked in the introduction, semisimple extrinsic symmetric triples are classified \cite{Nai}. Hence we may concentrate on extrinsic symmetric triples without simple ideals, which we will do in the second part of this section.

\subsection{Splitting into semisimple part and remainder}\label{Sssr}
\begin{lm} 
Each extrinsic symmetric triple $(\fg,\Phi,\ip)$ decomposes into a direct sum  $(\fg,\Phi,\ip)=(\fg_1,\Phi_1,\ip_1)\oplus (\fg_2,\Phi_2,\ip_2)$ of extrinsic symmetric triples, where $\fg_1$ is semisimple or zero and $\fg_2$ does not have simple ideals.
\end{lm}\label{Lsplit}
\proof Let $(\fg,\Phi,\ip)$  be an extrinsic symmetric triple. Since the sum of two semisimple ideals is semisimple $\fg$ has a unique semisimple ideal $\fs$. This must be invariant under automorphisms and derivations of $\fg$. In particular, $\fs$ is $\Phi$-invariant.  Let $\fr=\fs\cap\fs^\perp$ be the metric radical of $\ip|_{\fs\times \fs}$. Assume that $\fr\not=0$. Then $\fr$ is semisimple, hence $\fr=[\fr,\fr]$. Thus $\langle \fr,\fg\rangle=\langle[\fr,\fr],\fg\rangle = \langle\fr, [\fr,\fg]\rangle =0$ since $[\fr,\fg]\subset \fr$. This gives $\fr=0$. Thus $\fs$ is a nondegenerate extrinsic symmetric triple and we can decompose $\fg$ into $\fs\oplus\fs^\perp$. Since $\fs$ is maximal semisimple $\fs^\perp$ does not contain simple ideals. 
\qed

\subsection{Extrinsic symmetric triples without simple ideals}
In Section \ref{S3} we constructed $(\RR,\ZZ_2)$-equivariant metric Lie algebras $\dd$ by a certain extension procedure. One can prove that one can obtain in this way each $(\RR,\ZZ_2)$-equivariant metric Lie algebra $\fg$ that does not contain simple ideals (see \cite{KOesi}). Moreover, in $\fg=\dd$ the data $(\fl,\Phi_\fl)$ and $\fa$ can be choosen in such a way that the representation of $\fl$ on $\fa$ is semisimple. The proof is based on the existence of the so-called canonical isotropic ideal. The construction of this ideal goes back to an idea of L.\,B\'erard-Bergery. Let us recall briefly the definition and the properties of this ideal. For more detailed explanations see \cite{KO1, KOesi}.
Let $(\fg,\Phi,\ip)$ be an $(\RR,\ZZ_2)$-equivariant metric Lie algebra. There 
is a  chain of $\Phi$-invariant ideals 
$$\fg = R_{0}(\fg)\supset R_{1}(\fg)\supset R_{2}(\fg)\supset \ 
    \ldots \supset \ R_{l}(\fg)=0 $$  
which is defined by the condition that $R_{k}(\fg)$ is the smallest ideal of 
$\fg$ contained in $R_{k-1}(\fg)$ such that the $\fg$-module 
$R_{k-1}(\fg)/R_{k}(\fg)$ is semisimple.
The canonical ideal $\fri(\fg)\subset\fg$ is defined by 
$$ \fri(\fg):= \sum_{k=1}^{l} R_k(\fg) \cap  R_k(\fg)^{\perp} .$$ 
It has the following properties.
\begin{pr}{\rm \cite{KO1, KOesi}}
    If $(\fg,\Phi,\ip)$ is an $(\RR,\ZZ_2)$-equivariant  metric Lie 
    algebra, then $\fri(\fg)$ is a $\Phi$-invariant isotropic 
    ideal and the $\fg$-module $\fri(\fg)^{\perp}/\fri(\fg)$  is semisimple. 
    If $\fg$ does not contain simple ideals, then the Lie 
    algebra $\fri(\fg)^{\perp}/\fri(\fg)$ is abelian.
\end{pr}
We put  $\fl:=\fg/\fri(\fg)^\perp$ and $\fa:=\fri(\fg)^\perp/\fri(\fg)$. Let $\Phi_\fl$, $\Phi_\fa$ and $\ip_\fa$ be the structures induced by $\Phi$ and $\ip$, respectively. Then $\fri(\fg)$ can be identified with the abelian Lie algebra $\fl^*$. It follows that $(\fg,\Phi,\ip)$ is a quadratic extension of $(\fl,\Phi_\fl)$ by $\fa$ associated with a suitable $(\alpha,\gamma)\in\ZPhi$.

Similar to the situation of ordinary extensions of $\fl$ by $\fa$, which are classified by $H^2(\fl,\fa)$ equivalence classes of quadratic extensions of $(\fl,\Phi_\fl)$ by $\fa$ are in bijection to $\HPhi$. In general, an $(\RR,\ZZ_2)$-equivariant metric Lie algebra $(\fg,\Phi,\ip)$ without simple ideals can be isomorphic to different quadratic extensions $\fd_{\alpha,\gamma}(\fl,\Phi_{\fl},\fa)$ and $\fd_{\hat \alpha,\hat\gamma}(\hat\fl,\Phi_{\hat\fl},\hat\fa)$, where $\fl\not\cong\hat\fl$ or $\fa\not\cong\hat\fa$. This ambiguity vanishes if we only consider quadratic extensions $\dd$ for which the canonical isotropic ideal equals $\fl^*\subset \fl^*\oplus\fa\oplus\fl$. These extensions are called balanced. Cocycles $(\alpha,\gamma)\in\ZPhi$  for which $\dd$ is balanced are called balanced. In \cite{KOesi} we gave the following conditions  which allow to decide whether a cocycle is balanced or not. Note that the above considerations show that  the $(\fl,\Phi_\fl)$-module $\fa$ is semisimple if $\dd$ is balanced.
\begin{pr} Let $(\rho,\fa)$ be a semisimple orthogonal 
$(\fl,\Phi_{\fl})$-module and take $(\alpha,\gamma)\in \ZPhi$. 
Since $\rho$ is semisimple we have a decomposition $\fa=\fa^\fl\oplus 
\rho(\fl)\fa$ and a corresponding decomposition 
$\alpha=\alpha_0+\alpha_1$.  Let $m$ be such that $R_{m+1}(\fl)=0$. 
Then $(\alpha,\gamma)\in \ZPhi$ is  balanced if and only if it satisfies the 
following conditions $(A_{k})$ and $(B_{k})$ for $0\le k\le m$. 
\begin{enumerate} 
\item[$(A_0)$] 
Let $L_0\in \fz(\fl)\cap \ker \rho$ be such that there exist 
elements 
$A_0\in \fa$ and $Z_0\in \fl^*$ satisfying 
for all $L\in\fl$ 
\begin{enumerate} 
\item[(i)] $\alpha(L,L_0)=\rho(L) A_0 $, 
\item[(ii)] $\gamma(L,L_0,\cdot)=-\langle A_0,\alpha(L,\cdot)\rangle_\fa +\langle 
Z_0, [L,\cdot]_\fl\rangle$ as an element of $\fl^*$, 
\end{enumerate} 
then $L_0=0$. 
\item[$(B_0)$] The subspace $\alpha_0(\ker \lb_\fl)\subset 
\fa^\fl$ is non-degenerate. 
\item[$(A_k)$] $(k\ge 1)$\\ 
Let $\fk\subset S(\fl)\cap R_k(\fl)$ be an $\fl$-ideal such that there exist 
elements 
$\Phi_1\in \Hom(\fk,\fa)$ and $\Phi_2\in \Hom(\fk,R_k(\fl)^*)$ satisfying 
for all $L\in\fl$ and $K\in\fk$ 
\begin{enumerate} 
\item[(i)] $\alpha(L,K)=\rho(L)\Phi_1(K)-\Phi_1([L,K]_\fl)$, 
\item[(ii)] $\gamma(L,K,\cdot)=-\langle \Phi_1(K),\alpha(L,\cdot)\rangle_\fa 
+\langle 
\Phi_2(K), [L,\cdot]_\fl\rangle +\langle \Phi_2([L,K]_\fl), \cdot \rangle$ as an 
element of $R_k(\fl)^*$, 
\end{enumerate} 
then $\fk=0$. 
\item[$(B_k)$] $(k\ge 1)$\\ 
Let $\fb_k\subset\fa$ be the maximal submodule such that the system of equations 
$$ \langle\alpha(L,K), 
B\rangle_\fa=\langle\rho(L)\Phi(K)-\Phi([L,K]_\fl),B\rangle_\fa\ , 
\quad L\in\fl, K\in R_k(\fl), B\in\fb_k, $$ 
has a solution $\Phi \in \Hom(R_k(\fl),\fa)$. Then $\fb_k$ is non-degenerate. 
\end{enumerate} 
\end{pr} 

The property of a cocycle to be balanced only depends on its cohomology class. This yields the notion of a balanced cohomology class. We denote the set of balanced cohomology classes in $\HPhi$ by $\HPhi_b$. 

Let $(\fl,\Phi_{\fl})$ be an $(\RR,\ZZ_2)$-equivariant Lie algebra. We consider the 
category $M^{ss}_{\fl,\Phi_{\fl}}$ 
of semisimple orthogonal $(\fl,\Phi_{\fl})$-modules, where the morphisms between 
two modules 
$\fa_1,\fa_2$ are morphism of pairs $(S,U): ((\fl,\Phi_{\fl}),\fa_1)\rightarrow 
((\fl,\Phi_{\fl}),\fa_2)$. We denote the automorphism group of an object $\fa$ of 
$M^{ss}_{\fl,\Phi_{\fl}}$ 
by $G_{\fl,\Phi_{\fl},\fa}$. 
The natural right action of $G_{\fl,\Phi_{\fl},\fa}$ on $\HPhi$ leaves $\HPhi_{b}$ 
invariant.

Now we can formulate the following structure theorem.
 Let $\cL$ be a complete set of representatives of iso\-morphism classes of 
$(\RR,\ZZ_2)$-equivariant Lie algebras. 
For each $(\fl,\Phi_{\fl})\in\cL$ we choose a complete set of representatives 
$\cA_{\fl,\Phi_{\fl}}$ of isomorphism classes of objects in~$M^{ss}_{\fl,\Phi_{\fl}}$. 
\begin{theo} \label{Tstructure}
There is bijective map from  
 the union of orbit spaces 
 $$\coprod_{(\fl,\Phi_{\fl})\in\cL}\  \coprod_{\fa\in\cA_{\fl,\Phi_{\fl}}}\   
 \HPhi_{b}/ G_{\fl,\Phi_{\fl},\fa}\ .$$ 
to the set of isomorphism 
classes of $(\RR,\ZZ_2)$-equivariant metric Lie algebras without simple ideals.
This map sends the orbit of $[\alpha,\gamma]\in\HPhi_{b}$ 
to the isomorphism class of $\dd$. 
\end{theo}
We want to specialise this result to full extrinsic symmetric triples. We already noted that $\dd$ is $h$-graded if and only if $(\fl,\Phi_\fl)$ and $\fa$ are so. Furthermore, if $\fd:=\dd$ is balanced, then $[\fd^-,\fd^-]=\fd^+$ holds if and only if the conditions 
\begin{itemize}
\item[$(T_1)$] $\quad [\fl^-,\fl^-]=\fl^+$ and 
\item[$(T_{2})$] $\quad (\fa^\fl)^+=\alpha_0(\Ker \lb_{\fl^-})$
\end{itemize}
are satisfied, see \cite{KO2}, Lemma~5.5.~and Lemma~5.6. Here $\alpha_0$ denotes the orthogonal projection of $\alpha$ to $\fa^\fl$. 

Let $(\fl,\Phi_\fl)$ be an $h$-graded $(\RR,\ZZ_2)$-equivariant Lie algebra satisfying $[\fl^-,\fl^-]=\fl^+$ and let $\fa$ be an $h$-graded orthogonal  $(\fl,\Phi_\fl)$-module. We define 
$$\HPhi_\sharp:=\{[\alpha,\gamma]\in\HPhi_b\mid \alpha \mbox{ satisfies } (T_2)\}.$$ Then $\HPhi_\sharp$ is invariant under $G_{\fl,\Phi_{\fl},\fa}$. Let
$\cL^\sharp\subset \cL$ be the subset consisting of those $h$-graded $(\fl,\Phi_\fl)\in\cL$ with $[\fl^-,\fl^-]=\fl^+$. Furthermore, we denote by $\cA_{\fl,\Phi_{\fl}}^\sharp\subset\cA_{\fl,\Phi_{\fl}}$ the set consisting of those orthogonal  $(\fl,\Phi_\fl)$-modules in $\cA_{\fl,\Phi_{\fl}}$ that are $h$-graded. We conclude from Theorem~\ref{Tstructure}:
\begin{theo}\label{Tstructureii}
The map $\ZPhi \ni (\alpha,\gamma) \mapsto \dd$ induces a bijective map from 
 the union of orbit spaces 
 $$\coprod_{(\fl,\Phi_{\fl})\in\cL^\sharp}\  \coprod_{\fa\in\cA_{\fl,\Phi_{\fl}}^\sharp}\   
 \HPhi_{\sharp}/ G_{\fl,\Phi_{\fl},\fa}\ .$$ 
to the set of isomorphism 
classes of full extrinsic symmetric triples without simple ideals.
\end{theo}

Let us modify the cohomology set $\HPhi_\sharp$ used above to obtain an analogous result for indecomposable full extrinsic symmetric triples.
A non-trivial decomposition of a pair $((\fl,\Phi_{\fl}),\fa)$ 
consists of two non-zero morphisms of pairs 
$ (q_i,j_i):((\fl,\Phi_{\fl}),\fa)\longrightarrow 
((\fl_i,\Phi_{\fl_i}),\fa_i)$, $i=1,2,$ 
such that $(q_1\oplus q_2,j_1\oplus j_2):((\fl,\Phi_{\fl}),\fa)\rightarrow 
((\fl_1,\Phi_{\fl_1})\oplus (\fl_2,\Phi_{\fl_2}),\fa_1\oplus\fa_2)$ is 
an isomorphism. We will say that a cohomology class $\varphi\in 
\cH^2_Q(\fl,\Phi_{\fl},\fa)$ is decomposable if it can 
be written as a sum 
$$ \varphi=(q_1,j_1)^*\varphi_1 +(q_2,j_2)^*\varphi_2 $$ 
for a non-trivial decomposition $(q_i,j_i)$ of $((\fl,\Phi_\fl),\fa)$ and 
certain $\varphi_i\in \cH^2_Q(\fl_i,\Phi_{\fl_i},\fa_i)$, $i=1,2$. 
Here addition is induced by addition in the vector space $C^{2}(\fl,\fa)\oplus 
C^{3}(\fl)$.  
For $[\alpha,\gamma]\in\HPhi_\sharp$ the extrinsic symmetric triple $\dd$ is decomposable if and only if $[\alpha,\gamma]$ is decomposable \cite{KOesi}. Hence, if we replace $\HPhi_\sharp$ in Theorem~\ref{Tstructureii} by the set $\HPhi_0$ of all indecomposable elements of $\HPhi_\sharp$, then we get a bijection onto the set of isomorphism 
classes of indecomposable full extrinsic symmetric triples without simple ideals.

\subsection{The automorphism group}
Next we calculate the automorphism group of an extrinsic symmetric triple $(\fg,\Phi,\ip)$ in the case where $\fg=\dd$ is a balanced quadratic extension. 
\begin{pr}\label{PF}
Let $\dd$ be balanced. Then $F:\dd\rightarrow\dd$ is an automorphism of the extrinsic symmetric triple $\dd$ if and only if
\begin{equation}\label{EFF}
F=\left(\begin{array}{ccc} (S^{-1})^* & 0&0\\ 0&U^{-1} &0\\ 0&0&S \end{array}\right) \cdot\left(\begin{array}{ccc} \Id&-\tau^* & \bar\sigma-\frac12\tau^*\tau \\ 0&\Id&\tau\\ 0&0&\Id \end{array}\right):\fl^*\oplus\fa\oplus\fl\longrightarrow \fl^*\oplus\fa\oplus\fl
\end{equation}
and
\begin{itemize}
\item[(i)] $(S,U): ((\fl,\Phi_\fl),\fa)\rightarrow ((\fl,\Phi_\fl),\fa)$ is an isomorphism of pairs;
\item[(ii)] the bilinear map $\sigma:\fl\otimes\fl \rightarrow \RR$ defined by $\sigma(L_1,L_2):=\bar\sigma(L_1)(L_2)$ is antisymmetric and satisfies $(\tau,\sigma)\in{\cal C}^1_Q(\fl,\Phi_\fl,\fa)$; and
\item[(iii)] $(S,U)^*(\alpha,\gamma)=(\alpha,\gamma)(\tau,\sigma)^{-1}$.
\end{itemize}
\end{pr}
\proof Suppose that $F:\fd:=\dd\rightarrow \fd$ is an isomorphism. Since $\fd$ is balanced $F$ maps $\fri(\fd)=\fl^*$ to $\fri(\fd)=\fl^*$ and $\fri^\perp(\fd)=\fl^*\oplus \fa$ to $\fri^\perp(\fd)=\fl^*\oplus\fa$.  We define $S:=\proj_\fl\circ F|_\fl$ and $U:=(\proj_\fa\circ F|_\fa)^{-1}$. Then $(S,U): ((\fl,\Phi_{\fl}),\fa)   \rightarrow 
((\fl,\Phi_{\fl}),\fa)$ is a morphism of pairs and the map $F_1$ given by the first matrix in the proposition is an isomorphism from $\fd_{(S,U)^*(\alpha,\gamma)}(\fl,\Phi_{\fl},\fa)$ to $\fd$. Moreover, $F_1^{-1}\circ F:\fd \rightarrow \fd_{(S,U)^*(\alpha,\gamma)}(\fl,\Phi_{\fl},\fa)$ is an equivalence of quadratic extensions. Any such equivalence is given by a matrix that has the form of the second matrix in the proposition, where $\sigma$ and $\tau$ satisfy Conditions $(ii)$ and $(iii)$. These facts are explained in \cite{KOesi}. A detailed proof in the non-equivariant case can be found in \cite{KO1}.

Conversely, if $(S,U): ((\fl,\Phi_{\fl}),\fa)   \rightarrow 
((\fl,\Phi_{\fl}),\fa)$ is an isomorphism of pairs, then the map defined by the first matrix in the proposition is an isomorphism from $\fd_{(S,U)^*(\alpha,\gamma)}(\fl,\Phi_{\fl},\fa)$ to $\fd$. If $\sigma$ and $\tau$ satisfy Conditions $(ii)$ and $(iii)$, then the map defined by the second matrix is an equivalence from $\fd$ to $\fd_{(S,U)^*(\alpha,\gamma)}(\fl,\Phi_{\fl},\fa)$.
\qed

\begin{de}\label{DF}
Let be given an isomorphism $S:\fl\rightarrow \fl$, an isometry $U:\fa\rightarrow \fa$, a linear map $\tau:\fl\rightarrow \fa$, a cochain $\sigma\in C^2(\fl)$ and define $\bar\sigma:\fl\rightarrow \fl^*$ by $\bar\sigma(L_1)(L_2):=\sigma(L_1,L_2)$ for $L_1,L_2\in\fl$. Then we can define a linear map $F:\fl^*\oplus\fa\oplus\fl\rightarrow \fl^*\oplus\fa\oplus\fl$ by {\rm (\ref{EFF})}. We will denote this map by $F(S,U,\tau,\sigma)$.
\end{de}

\section{Classification of Lorentzian extrinsic symmetric spaces}
\subsection{Some known results} \label{Sknown}
Before classifying Lorentzian extrinsic symmetric spaces let us first recall some known results.  The case of parallel surfaces, i.e., of surfaces with parallel second fundamental form was considered by Chen, Dillen and Van der Veken. They classified parallel surfaces in 4-dimensional Lorentzian space forms  \cite{CV} and parallel Lorentzian surfaces in Lorentzian complex space forms \cite{CDV}.

Let us turn to the case of arbitrary dimension. If the associated extrinsic symmetric triple is semisimple, then we get a classification from Naitoh's list \cite{Nai}, see Section~\ref{Sss} for details. The main problem will be to classify extrinsic symmetric embeddings of non-semisimple symmetric spaces. Here as usual a symmetric space is called semisimple if its transvection group is semisimple. Let us recall the classification of (ordinary) non-semisimple Lorentzian symmetric spaces due to Cahen and Wallach \cite{CW}.   
Let   $(z,a_1,\dots,a_p,a_1',\dots,a_q',l)$ be the coordinates in $\RR\oplus\RR^{p+q}\oplus \RR$ and consider the metric 
$$g_{\lambda,\mu}:=2dzdl+ \sum_{i=1}^p da_i^2 + \sum_{j=1}^q {da'_j}^2 + \left(  \sum_{i=1}^p\lambda_i^2a_i^2 -  \sum_{j=1}^q\mu_j^2a_j'^2 \right)dl^2,$$
where 
$\lambda=(\lambda_1,\dots,\lambda_p)\in(\RR_{>0})^p$ and $\mu=(\mu_1,\dots,\mu_q)\in(\RR_{>0})^q$. Then
$M(p,q;\lambda,\mu):=(\RR^{p+q+2},\,g_{\lambda,\mu})$ is a Lorentzian symmetric space. In this notation we will omit $\lambda$ if $p=0$ and $\mu$ if $q=0$. Here these spaces are called Cahen-Wallach spaces. 
Conversely, every simply-connected indecomposable non-semisimple Lorentz\-ian symmetric space is either one-dimensional or isometric to one 
of the spaces $M(p,q;\lambda,\mu)$. Moreover, $M(p,q;\lambda,\mu)$ is isometric to $M(p',q';\lambda',\mu')$ if and only if $p=p'$, $q=q'$ and $(\lambda,\mu)=r(\lambda',\mu')$ for some $r>0$.

Every $n$-dimensional Cahen-Wallach space can be isometrically embedded into $\RR^{2,n}$. This was claimed in \cite{CW} and proven in detail by Blau, Figueroa-O'Farrill and Papadopoulos \cite{BFP}, who explicitly described the embedding as an intersection of two quadrics. It is easy to check that the embedding given in \cite{BFP} is extrinsic symmetric exactly for the Cahen-Wallach spaces $M(0,q;\mu)$ for $\mu=(1,\dots,1)$.

We want to end this section with explaining the relationship between 
full extrinsic symmetric triples and hermitian symmetric triples. Recall that a symmetric triple $(\fg,\theta,\ip)$ is the infinitesimal object that is associated with a pseudo-Riemannian symmetric space. It  consists of a metric Lie algebra $(\fg,\ip)$ and an isometric involution $\theta$ whose eigenspaces $\fg_+$ and $\fg_-$ satisfy $[\fg_-,\fg_-]=\fg_+$. The index of $\ip|_{\fg_-}$ is called index of the symmetric triple. It equals the index of the metric of any symmetric space associated with this triple. A hermitian symmetric triple $(\fg,D,\ip)$ is the infinitesimal object that is associated with a pseudo-hermitian symmetric space. It consists of a metric Lie algebra $(\fg,\ip)$ and an antisymmetric derivation $D$ satisfying $D^3=-D$ and a condition concerning the eigenspaces $\fg^+$ and $\fg^-$ of $\tau_D=\exp(\pi D)$. This condition is  $[\fg^-,\fg^-]=\fg^+$. In particular, $(\fg,\tau_D,\ip)$ is a symmetric triple and the almost complex structure $D|_{\fg^-}$ defines a hermitian structure on the associated symmetric space. We obtain the following relation between extrinsic and hermitian symmetric triples.
\begin{pr}\label{Phst}
If $(\fg,\Phi ,\ip)$, $\Phi=(D,\theta)$ is a full  extrinsic symmetric triple, then $(\fg,D,\ip)$ is a hermitian symmetric triple. \end{pr}

We will say that an extrinsic symmetric triple $(\fg,\Phi,\ip)$ is Lorentzian or of Lorentz type if the associated extrinsic symmetric space has Lorentz signature. If $(\fg,\Phi,\ip)$ is Lorentzian, then the hermitian symmetric triple $(\fg,D,\ip)$ has index two. 
Such hermitian triples were classified in \cite{KO2} and \cite{KOesi}. In order to determine all non-semisimple full extrinsic symmetric triples of Lorentz type we could first check the list in \cite{KO2} and determine all hermitian symmetric triples $(\fg,D,\ip)$ that admit an isometric involution $\theta$ anticommuting with $D$ and, second, classify these involutions. However, instead of checking lists we prefer to give a direct proof. 

Now let us start the classification of extrinsic symmetric triples and weak extrinsic symmetric triples of Lorentz type. Obviously, if such a triple is decomposable, then it decomposes into an indecomposable triple of Lorentz type and several indecomposable triples of Riemannian type, i.e., triples for which the associated extrinsic symmetric spaces are Riemannian ones. Thus to classify all full extrinsic symmetric triples of Lorentzian type we have to classify indecomposable ones of Lorentzian and Riemannian type. In both cases we have to distinguish between semisimple triples and those without simple ideals as shown in Section~\ref{Sssr}. This will be done in the following two sections. In Section~\ref{Sfwt} we will extend the classification to weak extrinsic symmetric triples.

\subsection{Classification of semisimple extrinsic symmetric triples}\label{Sss}

Semisimple extrinsic symmetric triples were classified by Naitoh \cite{Nai}. More exactly, Naitoh studied so-called semisimple symmetric graded Lie algebras, which are in a one-to-one correspondence with semisimple extrinsic symmetric triples. Let us briefly explain this correspondence. Roughly speaking, a semisimple symmetric graded Lie algebra is a semisimple graded Lie algebra $\bar \fg=\bar \fg_{-1}\oplus\bar \fg_0\oplus\bar \fg_{1}$ together with an involution $\bar \theta$ satisfying $\bar \theta(\bar \fg_\mu)=\bar \fg_{-\mu}$. The gradation of $\bar \fg$ defines a derivation $\bar D$, which anticommutes with $\bar \theta$. Let $\bar \fg=\bar \fg_+\oplus\bar \fg_-$ denote the decompositon of $\bar\fg$ defined by $\bar\theta$. The above mentioned correspondence  sends $\bar \fg$ to $\fg=\bar\fg_+\oplus i\bar\fg_-\subset\fg^{\Bbb C}$. Moreover, $\Phi=(D,\theta)$ is defined by $D=i\bar D$ and $\fg_+=\bar\fg_+$, $\fg_-=i\bar \fg_-$. 
On the other hand,  semisimple symmetric graded Lie algebras correspond to non-degenerate Jordan triple systems, which were classified by Neher \cite{Ne1, Ne2, Ne3}. Using this classification Naitoh gave a list of all semisimple symmetric graded Lie algebras. In particular we can read from this list all semisimple extrinsic symmetric triples of Lorentz type. We obtain four infinite series, for which we give $\fg,\fg_+$, $\fg^+$ and $\fg_+^+$ in the following table.   
\renewcommand{\arraystretch}{1.5}
$$ \mbox{\begin{tabular}{|c|c|c|c|}
\hline
$\fg$&$\fg_+$&$\fg^+$&$\fg_+^+$\\
\hline\hline
$\fsu(1,n)$&$\fso(1,n)$&$\frak u(1,n-1)$ &$\fso(1,n-1)$\\
\hline
$\fsu(2,n)$&$\fso(2,n)$&$\frak u(1,n)$ &$\fso(1,n)$\\
\hline
$\fso(3,k+l)$&$\fso(1,k)\oplus\fso(2,l)$&$\fso(2)\oplus\fso(1,k+l)$ &$\fso(0,k)\oplus\fso(1,l)$\\
\hline
$\fso(1,k+l+2)$&$\fso(0,k+1)\oplus\fso(1,l+1)$& $\fso(1,k+l)\oplus \fso(2)$&$\fso(0,k)\oplus\fso(1,l)$\\
\hline
\end{tabular}}$$

A classification of all non-compact  symmetric graded Lie algebras for which the involution is a Cartan involution was already achieved by  Kobayashi and Nagano \cite{KoNa}. Under the map $\bar\fg\mapsto \fg=\bar\fg_+\oplus i\bar\fg_-$ these symmetric graded Lie algebras correspond exactly 
to those extrinsic symmetric triples that are associated with embeddings of Riemannian symmetric spaces into Euclidean spaces. 

Semisimple extrinsic symmetric triples $(\hat \fg,\hat\Phi,\ip\hat{})$ that are associated with embeddings of Riemannian symmetric spaces into pseudo-Euclidean spaces can be obtained from those associated with embeddings into Euclidean spaces by the correspondence
$$\fg=\fg^+\oplus\fg^-\longmapsto\hat \fg=\fg^+\oplus i\fg^-\subset\fg^{\Bbb C}.$$ 
The involution $\hat \theta$ is given by the decomposition $\hat\fg=\hat \fg_+\oplus \hat \fg_-$, where $\hat \fg_+= \fg_+^+\oplus i\fg_+^-$ and $\hat\fg_-=\fg_-^+\oplus i\fg_-^-$. The derivation $D$ on $\fg$ is inner, i.e., $D=\ad(\xi)$. Then on $\hat\fg$ we also have $\hat D=\ad(\xi)$. In particular, the hermitian symmetric triple defined by $(\hat \fg,\ip\hat{},\hat\Phi)$ is the non-compact dual of the hermitian symmetric triple $\fg=\fg^+\oplus\fg^-$.

Each semisimple extrinsic symmetric triple is automatically full. Indeed, if $\fg$ is semisimple, then $\fg=[\fg,\fg]$. In particular, $\fg_-^+=[\fg_+^-,\fg_-^-]\oplus [\fg_+^+,\fg_-^+]$. Since $\fg_+^+=[\fg_+^-,\fg_+^-]$ the second summand is contained in $[[\fg_+^-,\fg_+^-],\fg_-^+]=[\fg_+^-,\fg_-^-]$. Thus  $\fg_-^+=[\fg_+^-,\fg_-^-]$.

\subsection{Classification of non-semisimple full extrinsic symmetric triples}\label{S5.3}
Now we consider extrinsic symmetric triples $(\fg,\Phi ,\ip)$ for which the Lie algebra $\fg$ does not have simple ideals. {}From Section~\ref{Sst} we know that every such triple  is a quadratic extension of an $(\RR,\ZZ_2)$-equivariant Lie algebra $(\fl,\Phi_\fl)$ by an orthogonal   $(\fl,\Phi_\fl)$-module $(\rho,\fa,\ip_\fa,\Phi_\fa)$. This gives immediately a classification in the case where the associated extrinsic symmetric space is full and Riemannian.

\begin{pr} \label{Pqr}
Any full Riemannian extrinsic symmetric triple $(\fg,\Phi ,\ip)$, $\Phi=(D,\theta)$, that does not have simple ideals is abelian. More exactly, $\fg$ is an abelian Lie algebra, $\ip$ is positive definite and $D$ satisfies $D^2=-\Id$. 
\end{pr}
\proof  As explained above $(\fg,\Phi ,\ip)$ is isomorphic to a quadratic extension $\fd:=\dd$. Since the associated extrinsic symmetric space is Riemannian we have $\fl^-=0$. Recall from Lemma~\ref{Lfest} that $\fg^+=[\fg^-,\fg^-]$ holds. This implies $\fl^+=[\fl^-,\fl^-]=0$. Hence $\fl=0$ and $\fd=\fa$. Moreover, $\fg^+=[\fg^-,\fg^-]$ implies $\fa=\fa^-$.
\qed

Let us turn to Lorentzian triples $(\fg,\Phi ,\ip)$. First we want to determine all $(\fl,\Phi_\fl)$ that can occur in the quadratic extension $\fg\cong\dd$. Let $\fh(1)$ denote the 3-dimensional Heisenberg algebra. This is a Lie algebra spanned by linearly independent vectors $X,Y,Z$ satisfying the relation $[X,Y]=Z$. We will use the notation 
$\fh(1)=\{[X,Y]=Z\}$. Similarly we write
\begin{eqnarray*}
\fsu(2)&=&\{[H,X]=2Y, [H,Y]=-2X, [X,Y]=2H\}, \\
\fsl(2,\RR)&=&\{[H,X]=2Y, [H,Y]=2X, [X,Y]=2H\}.
\end{eqnarray*}
For $\fl\in\{\fsu(2),\fsl(2,\RR)\}$ we denote by $\sigma_X, \sigma_Y, \sigma_H$ the basis of $\fl^*$ that is dual to to the basis $X,Y,H$ of $\fl$. Similarly, we write $\sigma_X, \sigma_Y, \sigma_Z$ for this basis if $\fl=\fh(1)$.
\begin{pr} \label{Pl}
Let $\fd:=\dd$ be a full Lorentzian extrinsic symmetric triple. Then either $\fl=0$ or $(\fl,\Phi_\fl)$ is isomorphic to one of the following $(\RR,\ZZ_2)$-equivariant Lie algebras: 
\begin{enumerate}
\item $\fl=\RR^2=\Span\{X,Y\}$,\\
$\Phi_\fl=(D_\fl,\theta_\fl)$ given by $D_\fl(X)=Y$, $D_\fl(Y)=-X$, $ \fl_+=\RR\cdot X ,\ \fl_-=\RR\cdot Y $;
\item $\fl=\fh(1)=\{[X,Y]=Z\}$,\\
$\Phi_\fl=(D_\fl,\theta_\fl)$ given by $\fl_+=\RR\cdot X$, $\fl_-=\Span\{Y,Z\}$,\\ $D_\fl(X)=Y$, $D_\fl(Y)=-X$, $D_\fl(Z)=0$;
\item $\fl=\fsu(2)=\{[H,X]=2Y,\ [H,Y]=-2X,\ [X,Y]=2H\}$,\\
$\Phi_\fl=(D_\fl,\theta_\fl)$ given by $\fl_+=\RR\cdot H$, $\fl_-=\Span\{X,Y\}$, $D_\fl=\frac12 \ad(X)$;
\item $\fl=\fsl(2,\RR)=\{[H,X]=2Y,\ [H,Y]=2X,\ [X,Y]=2H\}$, \\
$\Phi_\fl=(D_\fl,\theta_\fl)$ as in 3. 
\end{enumerate}
\end{pr}
\proof Since the associated extrinsic symmetric space has Lorentzian signature we have $\dim \fl_+^-=\dim \fl_-^-\le 1$. Now $[\fd^-,\fd^-]=\fd^+$ gives $\fl^+=[\fl^-,\fl^-]=[\fl_+^-,\fl_-^-]=\fl_-^+$.

If $\fl^-=0$ or $\fl^-$ is two-dimensional and abelian, this implies $\fl^+=0$. Hence $\fl=\fl^-$, which yields $\fl=0$ or $\fl\cong\RR^2$.

Now let us consider the case where $\fl^-$ is two-dimensional and non-abelian. We choose $L_+^-,L_-^-\in\fl$ such that 
$\fl_+^-=\RR\cdot L_+^-$, $\fl_-^-=\RR\cdot L_-^-$ and $L_-^-=D_\fl (L_+^-)$. Then $[L_+^-,L_-^-]=:L^+_-\not=0$ and $D_\fl (L_-^+)=0$. Moreover, we have $[L^+_-,L^-_+]=\lambda L_-^-$, $\lambda\in\RR$ and
$$[L^+_-,L^-_-]=[L^+_-,D_\fl( L^-_+)]=D_\fl [L^+_-,L^-_+]=\lambda D_\fl (L^-_-)=-\lambda L_+^-.$$
If $\lambda=0$ we put $X:=L^-_+$, $Y=L_-^-$, $Z=L^+_-$, from which we see that $\fl\cong\fh(1)$ and $D_\fl $ is given as in {\it 1}. If $\lambda\not=0$ we put
$$
X=(2/\lambda)\cdot L^+_-, \ \ Y=-(2/\textstyle{\sqrt {|\lambda|}})\cdot L^-_-, \ \ H=(2/\textstyle{\sqrt {|\lambda|}})\cdot L^-_+.
$$
For  $\lambda>0$ we obtain $\fl\cong\fsu(2)$ and $D_\fl $ as given  in {\it 3}. For $\lambda<0$ we have $\fl\cong\fsl(2,\RR)$.
\qed

The next question will be which semisimple representations of $(\fl,\Phi_\fl)$ can appear in the quadratic extension $\fg\cong\dd$ for those $(\fl,\Phi_\fl)$ determined in Prop.~\ref{Pl}.
\begin{pr}\label{Pn}
If $\fd:=\dd$ is a full Lorentzian extrinsic symmetric triple with $\fl\cong\RR^2$ or $\fl\cong\fh(1)$ and if $(\rho,\fa)$ is semisimple, then $\rho=0$.
\end{pr}
\proof If $\fl$ is solvable, then $\fd$ is also solvable. Now we use that $(\fd,D,\ip)$ is a hermitian symmetric triple by Prop.~\ref{Phst}. If the transvection group of a pseudo-hermitian symmetric space is solvable, then it is nilpotent \cite{KO2}, Prop.7.11. Hence, if $\rho$ is semisimple, then $\rho=0$. \qed

Now assume $\fl\in\{\fsu(2),\fsl(2,\RR)\}$. We put $\kappa=1$ if $\fl=\fsl(2,\RR)$ and  $\kappa=-1$ if $\fl=\fsu(2)$. We consider the adjoint representation on $\fa_3:= \fl$ and define $\Phi_\fa=(D_\fl,-\theta_\fl)$, $\hat \Phi_\fa=(D_\fl,\theta_\fl)$. Let $B_\fl$ be the Killing form on $\fl$ and put $\ip_\fa=\kappa B_\fl$. Then $\fa_3=(\ad,\fl,\ip_\fa,\Phi_\fa)$ and $\hat \fa_3=(\ad,\fl,\ip_\fa,\hat\Phi_\fa)$ are orthogonal $(\fl,\Phi_\fl)$-modules.

Furthermore, we consider the standard representation of $\fl$. For $\fl=\fsu(2)$ we have the standard representation $\rho_2$ on $\CC^2$.  Let $(\rho_4,\fa_4)$ be the realification of $(\rho_2,\CC^2)$. The real part of the standard scalar product on $\CC^2$ gives us a scalar product $\ip_\fa$ on $\fa_4$. For $\fl=\fsl(2,\RR)$ we have the standard representation $(\rho_2,\fv_2)$ on $\fv_2=\RR^2$. Using the dual representation $(\rho_2^*,\fv_2^*)$ we define $(\rho_4,\fa_4):=(\rho_2\oplus\rho_2^*,\fv_2\oplus\fv_2^*)$. 
In this case the dual pairing between $\fv_2$ and $\fv_2^*$ gives us an inner product $\ip_\fa$ on $\fa_4$. 

With respect to a suitable basis $a_1,\dots,a_4$ of $\fa_4$ the representation $\rho_4$ of $\fl$ is given by
$$\begin{array}{llll} 
X(a_1)=-a_3,& X(a_2)=-a_4,&X(a_3)=a_1,&X(a_4)=a_2,\\
Y(a_1)=a_4,& Y(a_2)=\kappa a_3,&Y(a_3)=a_2,&Y(a_4)=\kappa a_1,\\
H(a_1)=a_2,& H(a_2)=\kappa a_1,&H(a_3)=-a_4,&H(a_4)=-\kappa a_3.
\end{array}$$ 
and $\ip_\fa$ has the form $\diag(-\kappa,1,-\kappa,1)$. Moreover, we define $\Phi_\fa=(D_\fa,\theta_\fa)$ by 
$$\Ker D_\fa=\Span\{a_1,a_3\},\ D_\fa(a_2)=-a_4,\ (\fa_4)_+=\Span\{a_1,a_2\},\ (\fa_4)_-=\Span\{a_3,a_4\}.$$
Then $(\rho_4,\fa_4,\ip_\fa,\Phi_\fa)$ is an orthogonal $(\fl,\Phi_\fl)$-module.
 
 \begin{pr}\label{Psusl} Let $\dd$ be an indecomposable full Lorentzian extrinsic symmetric triple with $\fl\in\{\fsl(2,\RR),\ \fsu(2)\}$. Then, as an object of $M^{ss}_{\fl,\Phi_{\fl}}$, $\fa$ is isomorphic to $(\fa_3)^k\oplus(\hat \fa_3)^l\oplus(\fa_4)^m$  for exactly one triple  $(k,l,m)\in\NN^3$. 
 \end{pr}
\proof  We note that $\rho$ cannot have a trivial subrepresentation in $\fa$. Indeed, the maximal trivial subrepresentation $\fa^\fl$ is non-degenerate and invariant under $\Phi_\fa$. Moreover, $\fa^\fl$ is an ideal in $\dd$ since $[\alpha]\in H^2(\fa,\fl)=0$. Hence $\fa^\fl=0$ since otherwise $\dd$ would be decomposable.

Let $\fa_{\Bbb C}$ be the complexification of $\fa$. We extend $\ip_\fa$ to a hermitian inner product on $\fa _{\Bbb C}$.  In the following we will denote both a given linear map on $\fa$ and its complex linear extension to $\fa_{\Bbb C}$ by the same symbol.
We consider the map 
$$\begin{textstyle}B:=D_\fa-\frac12\rho(X)\end{textstyle}$$ on $\fa^{\Bbb C}$. By assumption $D_\fa$ is diagonalisable. Moreover, $\rho(X)$ is also diagonalisable since $X$ acts semi-simply on $\fl$. Both maps commute since $[D_\fa,\rho(X)]=\rho(D_\fl(X))=0$ by (\ref{EDrho}). Hence, as maps on $\fa_{\Bbb C}$ they have a common eigenspace decomposition. In particular, $B$ is diagonalisable. The eigenvalues of $B$ are purely imaginary since the eigenvalues of $D_\fa$ and $\rho(X)$ are so. 

Let $E_\lambda\subset\fa_{\Bbb C}$ denote the eigenspace with eigenvalue $\lambda$. The involution $\theta_\fa$ maps $E_\lambda$ to $E_{-\lambda}$ and complex conjugation also maps $E_\lambda$ to~$E_{-\lambda}$.

We now determine the eigenvalues of $B$. Each eigenspace $E_\lambda$ is an $\fl$-submodule since
$$\begin{textstyle} [B,\rho(L)]= [D_\fa-\frac12\rho(X),\rho(L)]=\rho(D_\fl(L))-\rho([\frac12 X,L])=\rho(D_\fl(L))-\rho(D_\fl(L))=0\end{textstyle}$$
for all $L\in\fl$. The maximal and the minimal eigenvalue of $\rho(X)|_{E_\lambda}$ have the form $\pm ki$ for some $k\in\NN$, where we can exclude $k=0$ because $\fa$ does not contain trivial subrepresentations. Since $D_\fa$ has only eigenvalues in $\{0,i,-i\}$ the condition  $B|_{E_\lambda}=\lambda\cdot\Id_{E_\lambda}$ yields
$$-(ki/2)+x_1=(ki/2)+x_2=\lambda$$
for suitable $x_1,x_2\in\{0,i,-i\}$. Hence $k\in\{1,2\}$. If $k=1$, then  $\lambda=\pm i/2$, if $k=2$, then $\lambda=0$.

Let us study the case $\fl=\fsu(2)$.  Take some $E_\lambda\not=0$. Then also $E_{-\lambda}\not=0$. By the above considerations $\lambda=\pm i/2$ or $\lambda=0$. If $\lambda=\pm i/2$, then $k=1$. Thus we can choose a submodule of $E_{i/2}$ that is isomorphic to the standard representation $(\rho_2,\CC^2)$. Let $v_+,v_-$ be a basis of this submodule such that 
$$X(v_\pm)=\mp v_\mp,\ \ Y(v_\pm)=i v_\mp,\ \
H(v_\pm)=\pm iv_\pm.$$ Then 
$(D_\fa-\frac12\rho(X)) v_\pm=\frac i2v_\pm$ implies
$$\begin{textstyle} D_\fa v_+=\frac i2 v_+-\frac 12 v_-, \  D_\fa v_-=\frac i2 v_-+\frac 12 v_+\,.\end{textstyle}$$
We define
\begin{eqnarray*}
v_1:=v_+ - v_-+\theta v_+-\theta v_-\,, && v_2:= \ i(v_++v_-+\theta v_++\theta v_-)\,,\\
v_3:= v_+ + v_--\theta v_+-\theta v_-\,,&&v_4:= -i(v_+-v_--\theta v_++\theta v_-)\,.
\end{eqnarray*}  
Since $v_1\not=0$ at least one of the vectors $\Re v_1$, $\Im v_1$ is not zero. If $\Re v_1\not=0$ we put $b_i:=\Re v_i$, otherwise $b_i:=\Im v_i$. Then the real vector space $W:=\Span\{b_1,\dots,b_4\}\subset \fa$ is invariant under $\rho$ and $\Phi_\fa$. Moreover, $W$ is non-degenerate since $\ip_\fa$ is positive definite on $\fa$. We may assume that $|b_1|=\dots=|b_4|=\sqrt2$.  If $U:W\rightarrow \fa_4$ is the linear map defined by $b_1\mapsto a_1+a_2$, $b_2\mapsto a_2-a_1$, $b_3\mapsto a_3+a_4$, $b_4\mapsto a_4-a_3$, then  $(\Id,U):((\fl,\Phi_\fl),\fa_4)\rightarrow ((\fl,\Phi_\fl),W)$  is an isomorphism of pairs. 

Now suppose that $\lambda=0$. Then $k=2$ and $\rho|_{E_\lambda}$ contains a multiple of the adjoint representation. Hence $\fn:=\Ker(\rho(X))\cap E_0\not=0$. In particular, at least one of the eigenspaces $\fn_+$, $\fn_-$ of $\theta|_\fn$ is not the zero space. Suppose $\fn_+\not=0$ and choose an element $a_X\in\fn_+$. Then $$W:=\Span\{a_X,\  a_Y:=(1/2)\cdot H(a_X),\  a_H:=-(1/2)\cdot Y(a_X)\}$$
is a three-dimensional subspace, which is invariant under $\rho$ and $\Phi$. Since $\ip_\fa$ is positive definite we may assume  $\langle a_X,a_X\rangle_\fa=8$. Let $U:W\rightarrow \fa_3$ be the linear map defined by $a_X\mapsto X$,  $a_Y\mapsto Y$,  $a_H\mapsto H$. Then $(\Id,U):((\fl,\Phi_\fl),\fa_3)\rightarrow ((\fl,\Phi_\fl),W)$  is an isomorphism of pairs. Analogously, for $a_X\in\fn_-$ we get an isomorphism  $(\Id,U):((\fl,\Phi_\fl),\hat \fa_3)\rightarrow ((\fl,\Phi_\fl),W)$. 

We proceed with $W^\perp$ instead of $\fa$, etc. The case $\fl=\fsl(2,\RR)$ is treated in the same way. In this case we use that $\ip_\fa$ is positive definite on $\fa^-$ to show that $\ip_\fa$ resticted to $W$ is non-degenerate.
\qed
\begin{theo}\label{Tfull}
If $(\fg,\Phi ,\ip)$, $\Phi=(D,\theta)$ is a non-semisimple indecomposable full extrinsic symmetric triple of Lorentz type, then it is isomorphic to exactly one quadratic extension $\dd$ associated with the data from the following list
\begin{enumerate}
\item \label{I1}$\fl=0$, $\fa=\RR^{2,0}$, $D_\fa(A_1)=A_2$, $D_\fa(A_2)=-A_1$,  $\fa_+=\Span\{A_1\}, \ \fa_-=\Span\{A_2\}$, where $A_1,A_2$ is the standard basis of $\RR^{2,0}$
\item \label{I2}$\fl=\RR^2=\Span\{X,Y\}$, $D_\fl(X)=Y$, $D_\fl(Y)=-X$, $ \fl_+=\RR\cdot X ,\ \fl_-=\RR\cdot Y, $
\begin{enumerate}
\item $\fa=\RR^{0,1}$, $D_\fa=0$, $\fa=\fa_-$, $\rho=0$,\\ $\alpha(X,Y)=A_0$, where $A_0$ is 
the standard basis of $\RR^{0,1}$,
$\gamma=0$;
\item $\fa=\RR^{1,0}$, $D_\fa=0$, $\fa=\fa_-$, $\rho=0$, \\$\alpha(X,Y)=A_0$, where $A_0$ is the standard basis of $\RR^{1,0}$,
$\gamma=0$;
\end{enumerate}
\item   \label{I3} $\fl=\fh(1)=\{[X,Y]=Z\},\ D_\fl(X)=Y$, $D_\fl(Y)=-X$, $D_\fl(Z)=0$\\
$ \fl_+=\RR\cdot X ,\ \fl_-=\Span\{Y,Z \}, $ 
\begin{enumerate}
\item[] $\fa=\RR^{2}$, $D_\fa(A_1)=A_2$, $D_\fa(A_2)=-A_1$,  $\fa_+=\RR\cdot A_2$, $\fa_-=\RR\cdot A_1$, where $A_1,A_2$ is the standard basis of $\RR^{2}$, $\rho=0$, \\$\alpha(X,Z)=A_1$, 
$\alpha(Y,Z)=A_2$,
$\alpha(X,Y)=0$, $\gamma=0$;
\end{enumerate}
\item   \label{I4}$\fl=\fsu(2)=\{[H,X]=2Y,\ [H,Y]=-2X,\ [X,Y]=2H\}$,  $D_\fl=\frac12 \ad(X)$,\\ $\fl_+= \RR\cdot H$, $\fl_-=\Span\{X, Y\}$,
\begin{enumerate}
\item[] $\fa=(\fa_3)^k\oplus(\hat\fa_3)^l\oplus(\fa_4)^m$, $\ k,l,m\in\NN$,\\
$\alpha=0$, $\gamma(H,X,Y)=4c,\ c\in\RR$; 
\end{enumerate}
\item \label{I5} $\fl=\fsl(2,\RR)=\{[H,X]=2Y,\ [H,Y]=2X,\ [X,Y]=2H\}$,  $\Phi_\fl=(D_\fl,\theta_\fl)$ as in 3.,
\begin{enumerate}
\item[] $\fa=(\fa_3)^k\oplus(\hat\fa_3)^l\oplus(\fa_4)^m$, $\ k,l,m\in\NN$,\\
$\alpha=0$,
$\gamma(H,X,Y)=4c,\ c\in\RR$.
\end{enumerate}
\end{enumerate}
\end{theo}
\proof We already know that $\fg\cong\dd$ for an  $(\RR,\ZZ_2)$-equivariant Lie algebra $(\fl,\Phi_\fl)$, an orthogonal $(\fl,\Phi_\fl)$-module $\fa$ and an indecompoable cohomology class $[\alpha,\gamma]\in\cH^2_Q(\fl,\Phi_\fl,\fa)_0/G_{\fl,\Phi_{\fl},\fa}$, where all these data are uniquely determined. In Propositions~\ref{Pl}, \ref{Pn} and \ref{Psusl} we determined all candidates for $(\fl,\Phi_\fl)$ and $\fa$. Let us compute $\cH^2_Q(\fl,\Phi_\fl,\fa)_0$ for these $(\fl,\Phi_\fl)$ and $\fa$. All these cohomologies can be easily computed directly. Here we will use our computations in \cite{KO1}.  
The case $\fl=0$ is trivial. Suppose $\fl=\RR^2$. By Prop.~\ref{Pn} the $\fl$-module $\fa$ is trivial. Since $\fg$ is indecomposable $\fa$ is in the image of $\alpha$,  thus either $\fa=0$ or $\fa$ is one-dimensional and $\fa=\fa^+_-$. Then \cite{KO1},~Lemma~6.8. implies that
$$C^2(\fl,\fa)\setminus \{0\} \ni \alpha \longmapsto [\alpha,0] \in\cH^2_Q(\fl,\Phi_\fl,\fa)_0$$
is a bijection. Hence $\fa\not=0$ and $\cH^2_Q(\fl,\Phi_\fl,\fa)_0/G_{\fl,\Phi_{\fl},\fa}$ consists of exactly one element.

Now take $\fl=\fh(1)$. Again $\fa$ is a trivial $\fl$-module by Prop.~\ref{Pn}.  Define 
$$Z_\fl:= \{\alpha \in C^2(\fl,\fa)\mid \alpha(X,Y) = 0,\, \alpha(X,Z)\in\fa_-^-,\, D(\alpha(X,Z))=\alpha(Y,Z)\}$$ 
By \cite{KO1}, Lemma~6.7.
$$\{\alpha \in Z_\fl \mid \alpha\not= 0,\ \alpha(\fl,\fl)=\fa \} \longrightarrow \cH^2_Q(\fl,\Phi,\fa)_0,\quad \alpha \longmapsto [\alpha,0]
$$
is a bijection. In particular, $\fa$ is two-dimensional and $\cH^2_Q(\fl,\Phi_\fl,\fa)_0/G_{\fl,\Phi_{\fl},\fa}$ consists of one element.

For $\fl\in \{\fsl (2,\RR), \fsu(2)\}$ we determined all possible $(\fl,\Phi_\fl)$-modules in Prop.~\ref{Psusl}. Let $\fa$ be such a module. Since $\fl$ is semisimple we have $H^2(\fl,\fa)$=0. Furthermore, $C^3(\fl)=\RR\cdot\gamma_0$ for $\gamma_0:=B_\fl(\lb,\cdot)$, where $B_\fl$ is the Killing form of $\fl$. Now one easily shows that 
$$C^3(\fl)\longrightarrow \cH^2_Q(\fl,\Phi,\fa)_0,\quad \gamma \longmapsto [0,\gamma]$$
is a bijection. Since in this case $G_{\fl,\Phi_{\fl},\fa}$ acts trivially on $\cH^2_Q(\fl,\Phi,\fa)_0$ the assertion follows.
\qed

If we integrate these triples, i.e., if we determine the associated extrinsic symmetric space $M_{\fg,\Phi}\hookrightarrow \fg_-$ for each triple $(\fg,\Phi,\ip)$ occuring in Theorem~\ref{Tfull}  we obtain the following classification result.
\begin{co} \label{Cfull} If $M\hookrightarrow V$ is a non-semisimple indecomposable full extrinsic symmetric space of Lorentz type, then it is isometric to exactly one of the following spaces:
\begin{enumerate}\label{Co}
\item \label{II1} $M=V=\RR^{1,0}$.
\item \label{II2} $V\in\{\RR^{1,2}, \RR^{2,1}\}$, where $\ip_V=2dx_1dx_3\pm dx_2^2$;\\
$M\cong\RR^{1,1}$ is the image of
$$\RR^2\longrightarrow V,\quad (r,s) \longmapsto (r,s^2,s).$$
\item \label{II3} $V= \RR^{2,3}$, where $\ip_V=2dx_1dx_4+2dx_2dx_5+ dx_3^2$;\\
$M\cong\RR^{1,2}$ is the image of
$$\RR^3\longrightarrow V,\quad (r,s,t) \longmapsto (s,-rt
+r^4/4,t,r,r^2).$$
\item \label{II4} $V= \RR^{2,2+k+2l+2m}$ with coordinates $(w_1,w_2,x,y,\hat y,z,\hat z,w_3,w_4)$, $x\in\RR^k,\ y,\hat y\in \RR^l,\ z,\hat z\in\RR^m$ and
$\ip_V=2dw_1dw_3+2dw_2dw_4+dx^2+dy^2+d\hat y^2+dz^2+d\hat z^2$; 

 $M$ is the image of 
\begin{eqnarray*}
\RR^{2+k+l+m}&\longrightarrow&  V \\
(r,t,s,v,u)&\longmapsto& \left(-({|s|^2}+|v|^2 +\textstyle{\frac12 }{|u|^2}+c)\cos r -(rc+t)\sin r+c,\right. \\
&& -({|s|^2}+|v|^2 +\textstyle{\frac12 }{|u|^2}+c)\sin r +(rc+t)\cos r,\ -s,\\
&& \left. v \cos r,\ -v \sin r,\ u \cos \textstyle{\frac r2},\  u\sin \textstyle{\frac r2},\ \frac12 (\cos r -1),\ \frac12\sin r \right),
\end{eqnarray*}
where $r,t\in\RR$, $s\in\RR^k$, $v\in\RR^l$, $u\in\RR^m$.

$M$ has signature $(1,1+k+l+m)$ and is isometric to a product $\RR^l\times \bar M(k,m)$, where $\bar M(k,m)$ is covered by the Cahen-Wallach space
$$M(0,k+m;(\underbrace{1,\dots,1}_{m\ \mbox{\footnotesize times}},\underbrace{2,\dots,2}_{k\ \mbox{\footnotesize times}})).$$
\item \label{II5} $V= \RR^{2+l+m,2+k+l+m}$ with coordinates $(w_1,w_2,x,y,\hat y,z,\hat z,w_3,w_4)$, $x\in\RR^k,\ y,\hat y\in \RR^l,\ z,\hat z\in\RR^m$ and
$\ip_V=2dw_1dw_3+2dw_2dw_4+dx^2+dy^2-d\hat y^2+dz^2-d\hat z^2$; 

$M$ is the image of 
\begin{eqnarray*}
\RR^{2+k+l+m}&\longrightarrow&  V \\
(r,t,s,v,u)&\longmapsto& \left(-({|s|^2}+|v|^2 +\textstyle{\frac12 }{|u|^2}-c)\cosh r -(rc+t)\sinh r-c,\right. \\
&& ({|s|^2}+|v|^2 +\textstyle{\frac12 }{|u|^2}-c)\sinh r +(rc+t)\cosh r,\ -s,\ v \cosh r,\\
&& \left.\ v \sinh r,\ u \cosh \textstyle{\frac r2},\  -u\sinh \textstyle{\frac r2},\ \frac12 (\cosh r -1),\ \frac12\sinh r \right),
\end{eqnarray*}
where $r,t\in\RR$, $s\in\RR^k$, $v\in\RR^l$, $u\in\RR^m$. 

$M$ has signature $(1,1+k+l+m)$ and is isometric to the product 
$$M(l,k,m):=\RR^l\times M(k+m,0; (\underbrace{1,\dots,1}_{m\ \mbox{\footnotesize times}},\underbrace{2,\dots,2}_{k\ \mbox{\footnotesize times}})).$$  

\end{enumerate}
\end{co}
\begin{re}{\rm
Note that in item \ref{II4} the Lorentzian manifold $\RR^l\times \bar M(0,m)$ (with its intrinsic metric) is a two-fold Riemannian covering of $\RR^l\times \bar M(m,0)$.
In item \ref{II5} the manifold $M(l,0,m)$ is even isometric to $M(l,m,0)$ (but not extrinsic isometric for $m\not=0$). Choosing  $k=m=0$ in items \ref{II4} and \ref{II5} we get indecomposable embeddings  of flat spaces $\RR^{1,1+l}\hookrightarrow \RR^{2,2+2l}$ and $\RR^{1,1+l}\hookrightarrow \RR^{2+l,2+l}$, respectively.}
\end{re}
{\sl Proof of Corollary \ref{Cfull}. } Recall from \cite{Kext1} that $M_{\fg,\Phi}\hookrightarrow \fg_-$ is the submanifold $G_+(0)$,  where $$G_+:=\langle\, \exp(\phi(X))\mid X\in\fg_+\,\rangle$$
with $\phi(X)=((\ad X)|_{\fg_-},-D(X))\in \fso(\fg_-)\ltimes\fg_- $.
The computation of items \ref{II1}, \ref{II2} and \ref{II3} from Thm.~\ref{Tfull}, \ref{I1} -- \ref{I3} is easy. Let us consider the cases \ref{II4} and \ref{II5}. Let $\tilde G_+$ be the simply connected Lie group  associated with $\fg_+$. This group is the transvection group of a Cahen-Wallach space. Cahen and Wallach \cite{CW} proved that $\tilde G_+$ can be identified with $\fg_+=\fl_+^*\oplus\fa_+\oplus\fl_+$ with group multiplication $$(Z,A,L)(\bar Z,\bar A,\bar L)=(Z+\bar Z+ \frac12 \left[e^{-\ad \bar L}(A),\bar A\right], e^{-\ad \bar L}(A)+\bar A,L+\bar L)$$ 
where $\lb$ and $\ad$ are the operations in $\fg_+$. See also \cite{Neu} for a detailed proof or \cite{KOesi} for a proof using the same notation as the present paper. Note that $\tilde G_+^+\cong \fa_+^+\subset \tilde G_+\cong \fl_+^*\oplus\fa_+\oplus\fl_+$ is the connected subgroup with Lie algebra $\fg_+^+\subset \fg$.  Moreover,   
$$s: \tilde G_+/\tilde G_+^+\longrightarrow \tilde G_+,\quad (Z,A^++A^-,L)\cdot \tilde G_+^+\longmapsto (Z+\textstyle{\frac12} [A^+,A^-], A^-,L),$$ with  $A^+\in\fa_+^+$, $A^-\in\fa_+^-$ is a global section of the projection $\tilde G_+\rightarrow \tilde G_+/\tilde G_+^+$. The image of $s$ equals $\fg_+^-\subset \tilde G_+$. Note that $(Z,A,L)=(0,0,L)(Z,A,0)=\exp(L)\exp(Z+A)$ holds in $\tilde G_+$. Hence $\fl_+^*\oplus\fa_+^-\oplus\fl_+\ni Z+A^-+L\mapsto \exp(L)\exp(Z+A)\in \tilde G_+$ is a bijection from the subspace $\fg_+^-$ of the Lie algebra $\fg_+$ onto the image of $s$. Since $\tilde G_+$ is a covering of $G_+$ we obtain  $M$ as the image of 
$$\iota: \fl_+^*\oplus\fa_+^-\oplus\fl_+\longrightarrow \fg_-,\quad  Z+A^-+L\longmapsto (\exp \Phi(L) \cdot \exp \Phi(Z+A))(0).$$ 
Now let us compute the image of $\iota$. Take an orthonormal basis $a_1,\dots,a_k$ of  $({\fa_3}_+^-)^k$,
an orthonormal basis $\hat a_1,\dots,\hat a_l$  of  $(\hat\fa_3{}_+^-)^l$ and an orthonormal basis $b_1,\dots,b_m$ of $(\fa_4{}_+^-)^m$. Then
$\sigma_X,\sigma_Y,\, \frac12X a_1,\dots,\frac12X a_k,\, \frac12X \hat a_1,\dots,\frac12X \hat a_l ,\, \frac\kappa 2Y \hat a_1,\dots,\frac\kappa 2Y \hat a_l ,-Xb_1,\dots,  -Xb_m,\linebreak \kappa Yb_1,\dots, \kappa Yb_m ,X,Y$
is an orthonormal basis of $\fg_-$. Using this basis we identify $(\fg_-,\ip|_{\fg_-})$ with $(V,\ip_V)$. We compute 
$$\Big(\exp \Phi(\frac r2 H)\cdot \exp \Phi(t\sigma_H +\sum_{i=1}^k s_i a_i +\sum_{i=1}^l v_i \hat a_i +\sum_{i=1}^m u_i b_i)\Big)(0)$$
with respect to the chosen basis and obtain the formulas claimed in items \ref{II4} and \ref{II5}.\qed

\begin{re}{\rm
A direct calculation shows that the mean curvature vector $h$ satisfies 
\begin{enumerate}
\item $h=0$ for all manifolds in Cor.~\ref{Cfull}, \ref{II1}--\ref{II3}, and 
\item $h=C\cdot \sigma_X$ with $C=-(4+2(k+l)+m))/(2+k+l+m)$ 
for all manifolds in Cor.~\ref{Cfull}, \ref{II4}--\ref{II5}, including the flat ones.
\end{enumerate}}
\end{re}
\subsection{Extensions to full weak triples} \label{Sfwt}
\subsubsection{The problem}\label{Sfwt1}

In Section \ref{S5.3} we classified all Lorentzian extrinsic symmetric spaces that are full in some pseudo-Euclidean space. Now we want to extend this classification to the case of Lorentzian extrinsic symmetric spaces $M\hookrightarrow\RR^{p,q}$ for which the minimal affine subspace $V$ containing $M$ can be degenerate. As explained in \cite{Kext1} it suffices to classify full Lorentzian extrinsic symmetric spaces in (possibly degenerate) inner product spaces. This classification problem is equivalent to the classification of full weak extrinsic symmetric triples of Lorentz type. In the following we will often denote a weak extrinsic symmetric triple just by the symbol of the underlying Lie algebra.
Let us recall the relation of extrinsic symmetric triples and weak extrinsic symmetric triples from \cite{Kext1}.
 
Let $R$ be a finite-dimensional real vector space. We consider $R$ as an abelian Lie algebra. Then $R$ together with the trivial inner product $\ip=0$, the trivial derivation $D=0$ and the involution $\theta=-\Id$ can be understood as a weak extrinsic symmetric triple. A {\it normal extension} of an extrinsic symmetric triple $(\fg,\Phi,\ip)$ by $R$ is a short exact sequence of weak extrinsic symmetric triples
$$0\longrightarrow R \longrightarrow \tilde\fg\longrightarrow \fg\longrightarrow 0.$$
Each such sequence is a central extension of the underlying Lie algebras. In particular, equivalence classes of normal extensions of $(\fg,\Phi,\ip)$ by $R$ are classified by $H^2(\fg,R)^D_-:=\{a\in H^2(\fg,R)\mid \theta(a)=-a,\ Da=0\}$, where $R$ is considered as a trivial $\fg$-module. The weak extrinsic symmetric triple $\tilde \fg$ is full  if and only if $\fg$ is full and if the associated cohomology class $[\omega]$ satisfies $\omega(\Ker \lb|_{\fg^-_-\wedge\fg^-_+})=R$.

If $(\tilde\fg,\tilde\Phi,\ip\tilde{}\,)$ is a weak extrinsic symmetric triple and $R:=\tilde\fg^\perp$ is its metric radical, then  $\tilde\fg/R$ together with the structure $(\Phi,\ip)$ induced by $(\tilde\Phi,\ip\tilde{}\,)$ is an extrinsic symmetric triple and  $(\tilde\fg,\tilde\Phi,,\ip\tilde{}\,)$ is a normal extension of $(\fg,\Phi,\ip)$  by $R$. In particular, isomorphism classes of weak extrinsic symmetric triples $(\tilde\fg,\tilde\Phi,\ip\tilde{}\,)$ with metric radical $R$ and $\tilde\fg/R\cong \fg$ (as extrinsic symmetric triples) correspond bijectively to the elements of 
 $H^2(\fg,R)^D_-/(\Aut(\fg,\Phi,\ip)\times \GL( R))$.  
 
The latter result together with Lemma \ref{Lsplit} shows that each weak extrinsic symmetric triple $(\fg,\Phi,\ip)$ decomposes into a direct sum  $(\fg,\ip,\Phi)=(\fg_1,\Phi_1,\ip_1)\oplus (\fg_2,\Phi_2,\ip_2)$, where $(\fg_1,\Phi_1,\ip_1)$ is a semisimple extrinsic symmetric triple or zero and $(\fg_2,\Phi_2,\ip_2)$ is a weak extrinsic symmetric triple which arises as a central extension of an extrinsic symmetric triple without simple ideals.

 \subsubsection{Second cohomology and derivations}

We denote by $\Der(\fg)$ the Lie algebra of antisymmetric derivations of the metric Lie algebra $(\fg,\ip)$. Moreover, we consider
$$\Der(\fg)^D_-:=\{ \ph\in \Der(\fg)\mid \ph\circ D=D\circ \ph,\ \theta\circ \ph=-\ph\circ\theta\}$$
and the Lie algebra of outer derivations $$\Out(\fg)^D_-:=\Der(\fg)^D_-/\ad(\fg_-^+).$$
\begin{lm} \label{LDer}
There is an isomorphism
$H^2(\fg, \RR)^D_-\cong\Out(\fg)^D_-$.
 \end{lm}
\proof We consider the map
\begin{equation}\label{Eio}
\iota:\ \Der(\fg)\ni\ph \longmapsto \omega_\ph:=\langle\ph(\cdot),\cdot\,\rangle\in Z^2(\fg,\RR).
\end{equation}
Since $\ip$ is non-degenerate this map is an isomorphism. Moreover,
\begin{eqnarray*}
\omega \in B^2(\fg,\RR)& \Leftrightarrow &\exists\,\eta\in C^1(\fg,\RR) :\omega=d\eta\ \ \Leftrightarrow 
\ \ \exists\,\eta\in C^1(\fg,\RR) :\omega=-\eta(\lb)\\
& \Leftrightarrow &\exists\,X\in\fg :\omega=\langle X,\lb\rangle \ \ \Leftrightarrow 
\ \ \exists\,X\in\fg :\omega=\langle \ad(X)(\cdot),\cdot\,\rangle\\
& \Leftrightarrow & \omega \in \iota(\ad(\fg))
\end{eqnarray*}
The conditions $D\circ \ph=\ph\circ D$ and $\theta\circ \ph=-\ph\circ\theta$ translate into $D\omega_\ph=0$ and $\theta^*\omega_\ph=-\omega_\ph$,
which completes the proof of the proposition. \qed

\begin{pr} \label{Pphi} Let $\fg=\dd$ be balanced. 
Then $\ph:\fg\rightarrow\fg$ is an element of $\Der(\fg)^D_-$ if and only if
\begin{equation}\label{Eph}
\ph= \left(\begin{array}{ccc} -\hat S^*&-\hat\tau^* & \hat\sigma \\ 0&-\hat U&\hat \tau\\ 0&0&\hat S \end{array}\right):\fl^*\oplus\fa\oplus\fl\longrightarrow \fl^*\oplus\fa\oplus\fl
\end{equation}
and
\begin{itemize}
\item[(i)] $\hat S:\fl\rightarrow \fl$ is a derivation satisfying $D_\fl\circ\hat S=\hat S\circ D_\fl$ and $\theta_\fl\circ\hat S=- \hat S\circ \theta_\fl$,\\
$\hat U:\fa\rightarrow \fa$ is a linear antisymmetric map that satisfies $D_\fa\circ\hat U=\hat U\circ D_\fa$ and $\theta_\fa\circ\hat U=- \hat U\circ\theta_\fa$, and, moreover
 $$ \rho(\hat S(L))=[\rho(L),\hat U] $$
holds for all $L\in\fl$;
\item[(ii)] the bilinear map $\hat\sigma:\fl\otimes\fl \rightarrow \RR$ defined by $\hat\sigma(L_1,L_2):=\hat\sigma(L_1)(L_2)$ is antisymmetric and $\hat\tau\in C^1(\fl,\fa)$, $\hat \sigma\in C^2(\fl,\fa)$ satisfy
$$
(\theta_\fl,\theta_\fa)^*\hat \tau=-\hat \tau,\ \ \theta_\fl^*\hat\sigma=-\hat\sigma,\ \
(e^{tD_\fl},e^{-tD_\fa})^*\hat \tau=\hat \tau,\ \ (e^{tD_\fl})^*\hat\sigma=\hat\sigma,
$$
for all $t\in\RR$; and
\item[(iii)] the equations \\[1ex]
$\hat U(\alpha(L_1,L_2))+\alpha(\hat S L_1,L_2)+\alpha(L_1,\hat S L_2)=-d\hat \tau (L_1,L_2)$\\[0.5ex]
$\gamma(\hat S L_1,L_2,L_3)+\gamma(L_1,\hat S L_2,L_3)+\gamma(L_1,L_2,\hat S L_3)= -(d\hat\sigma +\langle\alpha\wedge \hat\tau\rangle)(L_1,L_2,L_3)$\\[1ex]
hold for all $L_1,L_2,L_3\in\fl$.
\end{itemize}
\end{pr} 
\proof  Roughly speaking, $\Der(\fg)^D_-$ differs from the Lie algebra of the automorphism group $\Aut(\fg)$ of the extrinsic symmetric triple $\fg$ just by the commutation rule for $\theta$. Hence we can determine $\Der(\fg)^D_-$ by differentiating Equation (\ref{EFF}) and the Conditions (i), (ii), (iii) in Prop.~\ref{PF} and by changing the sign in all corresponding commutation rules for $\theta_\fl$ and $\theta_\fa$.
Of course, the assertion can also be easily verified by a direct computation. \qed
\begin{de}\label{Dph}
Let be given a  derivation $\hat S:\fl\rightarrow \fl$, an antisymmetric linear map $\hat U:\fa\rightarrow \fa$, a linear map $\hat \tau:\fl\rightarrow \fa$, a cochain $\hat \sigma\in C^2(\fl)$ and define $\hat\sigma:\fl\rightarrow \fl^*$ by $\hat\sigma(L_1)(L_2):=\hat\sigma (L_1,L_2)$ for $L_1,L_2\in\fl$. Then we can define a linear map $\ph:\fl^*\oplus\fa\oplus\fl\rightarrow \fl^*\oplus\fa\oplus\fl$ by {\rm (\ref{Eph})}. We will denote this map by $\ph(\hat S,\hat U,\hat \tau,\hat \sigma)$.
\end{de}
The action of  $\Aut(\fg)$ on $Z^2(\fg,\RR)$ translates into the adjoint action on $\Der(\fg)$. More exactly, the map $\iota$ defined in (\ref{Eio}) has the property 
$$F^*(\iota(\ph))=\iota(F^{-1}\ph F)$$
for all $F\in \Aut(\fg)$ and $\ph\in\Der(\fg)$.

\subsubsection{Decomposability}
A weak extrinsic symmetric triple is called decomposable if it is the direct sum of two non-trivial weak extrinsic symmetric triples. Otherwise it is called indecomposable. Any weak extrinsic symmetric triple of Lorentz type decomposes into an indecomposable weak extrinsic symmetric triple of Lorentz type and a (possibly trivial) Riemannian weak extrinsic symmetric triple. Obviously, the triple is full if and only if both summands are full. 
Since our aim is the classification of all full weak extrinsic symmetric triples of Lorentzian type we have to classify all indecomposable and full weak extrinsic symmetric triples which have Lorentzian or Riemannian type.  

Let $\tilde \fg$ be an indecomposable and full weak extrinsic symmetric triple of Lorentzian or Riemannian type. As above we define $R:=\tilde\fg^\perp$. Then the extrinsic symmetric triple $\tilde \fg/R$ is not necessarily indecomposable. However, the following holds.

\begin{pr}\label{Prt}
A weak extrinsic symmetric triple $(\tilde\fg,\tilde\Phi,\ip\tilde{}\,)$ of Riemannian type is full and indecomposable if and only if 
\begin{enumerate}
\item  $\tilde\fg$ is an indecomposable semisimple extrinsic symmetric triple; or 
\item $\fg :=\tilde \fg/R = \fa_0=\fa_0^-$ is abelian and and has a positive definite scalar product; and
the cohomology class $[\omega]\in H^2(\fg,R)^D_-$ defined by the sequence $0\rightarrow R \rightarrow \tilde\fg\rightarrow \fg\rightarrow 0$ does not decompose as a sum $[\omega']+[\omega'']$ for which
$$
\omega'(\fa_0',\fa_0')\subset R', \quad \omega'(\fa_0'',\fa_0)=0,\quad 
\omega''(\fa_0'',\fa_0'')\subset R'', \quad \omega''(\fa_0',\fa_0)=0
$$
for some decompositions $\fa_0=\fa_0'\oplus\fa_0''$, $R=R'\oplus R''$ with $\fa_0'\oplus R'\not=0$ and $\fa_0''\oplus R''\not=0$.\qed\end{enumerate}
\end{pr}

\begin{pr}\label{Plt}  A weak extrinsic symmetric triple $(\tilde\fg,\tilde\Phi,\ip\tilde{}\,)$ of Lorentz type is full and indecomposable if and only if \begin{enumerate}
\item  $\fg:=\tilde\fg/R=\fg_1\oplus \fa_0$ (direct sum of extrinsic symmetric triples), where $\fg_1$ is indecomposable of Lorentz type and  $\fa_0=\fa_0^-$ is abelian and positive definite; and
\item the cohomology class $[\omega]\in H^2(\fg,R)^D_-$ defined by the sequence $0\rightarrow R \rightarrow \tilde\fg\rightarrow \fg\rightarrow 0$ is not in the $\Aut(\fg_1\oplus\fa_0)$-orbit of a class $[\omega']+[\omega'']$ satisfying
$$
\omega'(\fg_1\oplus\fa_0',\fg_1\oplus\fa_0')\subset R', \quad \omega'(\fa_0'',\fg)=0,\quad 
\omega''(\fa_0'',\fa_0'')\subset R'', \quad \omega''(\fg_1\oplus \fa_0',\fg)=0
$$
for some decompositions $\fa_0=\fa_0'\oplus\fa_0''$, $R=R'\oplus R''$ with $\fa_0''\oplus R''\not=0$; and
\item $\fg_1$ is full.
\end{enumerate}
\end{pr}
\proof
Let $(\tilde\fg,\tilde\Phi,\ip\tilde{}\,)$ be full and indecomposable. As explained in Section \ref{Sknown} we have $\tilde\fg/R=\fg_1\oplus\fg_2$, where $\fg_1$ is a full and indecomposable extrinsic symmetric triple of Lorentzian type and $\fg_2$ is a full extrinsic symmetric triple of Riemannian type. The Riemannian part $\fg_2$ cannot contain simple ideals. Indeed, otherwise $\fg_2=\fg_2'\oplus \fg_2''$, where $\fg_2''$ is semisimple. In particular, $H^1(\fg_2'',R)=H^2(\fg_2'',R)=0$. This would imply that $\tilde \fg$ is decomposable since under these assumptions  $H^2(\tilde\fg/R,R)=H^2(\fg_1\oplus\fg_2',R)$ by the K\"unneth formula. Hence  $\tilde\fg$ would be decomposable into $\fg''_2$ and an extension of $\fg_1\oplus\fg_2'$. Thus $\fg_2$ does not contain simple ideals, hence $\fg_2=:\fa_0$  is abelian by Prop.~\ref{Pqr}.  If $[\omega]$ were in the $\Aut(\fg_1\oplus\fa_0)$-orbit of a class $[\omega']+[\omega'']$ as in {\it 2.}, there would exist a morphism of short exact sequences 
\begin{equation}\label{Eseq}
\begin{array}{ccccccccl} 0&\rightarrow &R&\rightarrow &\tilde\fg&\rightarrow& \fg_1\oplus\fa_0&\rightarrow&0\\[0.5ex]
& &\downarrow& &\downarrow&&\downarrow&&\\[0.5ex]
0&\rightarrow &R'\oplus R''&\rightarrow &\tilde\fg'\oplus\tilde\fg''&\rightarrow& (\fg_1\oplus\fa_0')\oplus \fa_0''&\rightarrow&0\,,
\end{array}
\end{equation}
where all vertical maps are isomorphisms. Hence $\tilde\fg$ would be decomposable.

Now let $\tilde\fg$ be decomposable or not full and assume that the conditions {\it 1.} and {\it 3.} hold. We have to show that condition {\it 2.} cannot be satisfied. If $\tilde \fg$ is decomposable, then $\tilde \fg= \tilde \fg'\oplus \tilde\fg''$, where $\tilde \fg'$ is a weak extrinsic symmetric triple of Lorentzian type and  $\tilde \fg''$ is a weak extrinsic symmetric triple of Riemannian type. Let $R=\tilde\fg\cap\tilde\fg^\perp$,   $R'=\tilde\fg'\cap(\tilde\fg')^\perp$ and $R''=\tilde\fg''\cap(\tilde\fg'')^\perp$ be the metric radicals of $\tilde\fg$, $\tilde \fg'$ and $\tilde\fg''$, respectively. Then $R=R'\oplus R''$ and $\fg_1\oplus\fa_0=\tilde\fg'/R'\oplus \tilde\fg''/R''$. Since the summands of any decomposition of $\fg_1\oplus\fa_0$ into indecomposable extrinsic symmetric triples are unique up to isomorphisms we get $\tilde\fg'/R'\cong \fg_1\oplus\fa_0'$ and $\tilde\fg''/R''\cong\fa_0''$, where $\fa_0=\fa_0'\oplus\fa''_0$. Hence there exists a morphism of short exact sequences as in (\ref{Eseq}). In particular, $[\omega]$ is in the $\Aut(\fg_1\oplus\fa_0)$-orbit of a class $[\omega']+[\omega'']$, where 
$\omega'(\fg_1\oplus\fa_0',\fg_1\oplus\fa_0')\subset R'$, $\omega'(\fa_0'',\fg)=0$, 
$\omega''(\fa_0'',\fa_0'')\subset R''$, $\omega''(\fg_1\oplus \fa_0',\fg)=0$. Now suppose that $\tilde \fg$ is not full and assume that the conditions {\it 1.} and {\it 3.} hold. We will show that in this case $\tilde \fg$ is decomposable, which will finish the proof. If $\tilde \fg$ is not full, then there exists an element $X\in\tilde\fg_-^+$ such that $X\not\in[\tilde\fg_+^-,\tilde\fg_-^-]$. Since $\tilde\fg/R$ is full there is an element $r\in R$ such that $X+r 
\in[\tilde\fg_+^-,\tilde\fg_-^-]$ thus $r\not \in[\tilde\fg_+^-,\tilde\fg_-^-]$. Now choose a subspace $U\subset\tilde\fg$ that is complementary to $\RR \cdot r$ in $\tilde \fg^+_-$ and that contains $[\tilde\fg_+^-,\tilde\fg_-^-]$. Then $\fg_2:=\tilde\fg_+\oplus U\oplus\tilde\fg_-^-$ is a weak extrinsic symmetric triple. Hence $\tilde\fg= \fg_2\oplus \RR\cdot r$ is a direct sum of non-trivial weak extrinsic symmetric triples.
\qed

Next we determine $H^2(\fg_1\oplus\fa_0,R)$ for all relevant cases in Prop.~\ref{Prt} and Prop.~\ref{Plt}. 
\subsubsection{Riemannian type} \label{SSrt}
As in the previous section let $\fa_0=\fa_0^-$ be abelian and let the scalar product on $\fa_0$ be positive definite. Let $S^2(V,W)$ denote the space of symmetric bilinear maps from $V\times V$ to~$W$.
\begin{pr}
We have $H^2(\fa_0,R)^D_- \cong S^2((\fa_0)_-,R)$.
\end{pr}
\proof Obviously,
$\Out(\fa_0)^D_-=\fso(\fa_0)^D_-:=\{\ph\in\fso(\fa_0)\mid \ph\circ \theta=-\theta\circ \ph,\ \ph\circ D=D\circ\ph\}$
and
\begin{equation}\label{ESym}
\fso(\fa_0)^D_-\longrightarrow \Sym^2((\fa_0)_-^*),\quad \ph\longmapsto \langle \ph\circ D(\cdot),\cdot\,\rangle \end{equation}
 is an isomorphism. This gives $H^2(\fa_0,R)^D_- \cong \Out(\fa_0)^D_-\otimes R\cong \Sym^2((\fa_0)_-^*)\otimes R\cong S^2((\fa_0)_-,R)$.
\qed

This immediately implies:
\begin{pr}  \label{PRt} There is a bijection from
$H^2(\fa_0,R)_-^D/(\Aut(\fa_0)\times \GL(R))$ to
$$S^2((\fa_0)_-,R)/(\grO((\fa_0)_-)\times\GL(R)),$$
where $\grO((\fa_0)_-)\times\GL(R)$ acts in the natural way. 
\end{pr}
\begin{de} An element $B\in S^2((\fa_0)_-,R)$ is called decomposable if there are decompositions $R= R'\oplus R''$ and $(\fa_0)_-=\fa'\oplus\fa''$, $\fa'\perp\fa''$ with $\fa'\oplus R'\not=0$ and $\fa''\oplus R''\not=0$ such that  
$$B(\fa',\fa')\subset  R',\ B(\fa'',\fa'')\subset R'',\ B(\fa',\fa'')=0.$$
\end{de}
Obviously, decomposability of $B$ only depends on the $\Aut(\fa_0)\times \GL(R)$-orbit $[B]$ of~$B$, hence we can speak of decomposable and indecomposable orbits. Prop.~\ref{Prt} now gives:
\begin{pr}
\quad The weak extrinsic symmetric triple corresponding to the orbit $[B]$ is full and indecomposable if and only if $[B]$ is indecomposable.
\end{pr}
\subsubsection{Lorentz type}
We have to determine $H^2(\fg,R)^D_-$ for all $\fg=\fg_1\oplus \fa_0$ with $\fg_1=\dd$ as listed in Thm.~\ref{Tfull}. Since $\dd$ is balanced we can understand $\fg$ as a balanced quadratic extension $\fd_{\alpha,\gamma}(\fl,\Phi_\fl,\fa\oplus\fa_0)$, where $\fa_0$ is considered as trivial $(\fl,\Phi_\fl)$-module.
\paragraph{The case \fbox{$\fl=0$}} 
Here $\fg$ is isomorphic to the abelian metric Lie algebra $\fa^{2,0}:=\RR^{2,0}$ with $D_\fa$ and $\theta_\fa$ determined by $\fa^{2,0}=(\fa^{2,0})^-$. Thus the problem reduces to determine $H^2(\fa_0)^D_-$ for an abelian Lie algebra $\fa_0=\fa_0^-$ with scalar product of index two. The result is exactly the same as in Section \ref{SSrt} including indecomposability.

\paragraph{The case \fbox{$\fl=\RR^2$}}

Let $\fg_1:=\dd$ be as defined in Thm.~\ref{Tfull}, \ref{I2}.~(a) or (b). 
\begin{lm}\label{LHR2} For $\fl=\RR^2$ and $\fg=\fg_1\oplus\fa_0$ with $\fg_1:=\dd$ as considered above 
$$H^2(\fg,R)_-^D\,\cong\, R\,\oplus \,S^2((\fa_0)_-,R)\,\oplus\, (\fa_0)_-\otimes R.$$
\end{lm}
\proof We first compute $\ad(\fg_-^+)$. Let $\sigma_X, \sigma_Y$ be the basis of $\fl^*$ that is dual to the basis $X,Y$ of $\fl$. We have $\fg^+_-=\RR\cdot A_0$ and 
$$[A_0,\fl^*\oplus\fa\oplus\fa_0]=0,\ [A_0,X]=\pm \sigma_Y,\ [A_0,Y]=\mp \sigma_X,$$
where the sign depends on the choice of the scalar product on $\fa$ according to Thm.~\ref{Tfull}, \ref{I2}.~(a) or (b). Hence
$$\ad(\fg^+_-)=\{\ph\in\Der(\fg)^D_-\mid \ph=\ph(0,0,0,\hat\sigma)\}$$
and, consequently,
$$\Out(\fg)^D_-\cong\{\ph(\hat S,\hat U,\hat \tau,\hat \sigma)\in\Der(\fg)^D_-\mid \hat \sigma=0\}.$$
Let us check which $\ph(\hat S,\hat U,\hat \tau,\hat 0)$ are in $\Der(\fg)^D_-$. The linear map $\hat S:\fl\rightarrow\fl$ is a derivation satisfying $D_\fl\circ\hat S=\hat S\circ D_\fl$ and $\theta_\fl\circ\hat S=- \hat S\circ \theta_\fl$ if and only if $\hat S(X)=\mu Y$ and $\hat S(Y) =-\mu X$ for some $\mu\in\RR$. 
Furthermore, $\hat U:\fa\oplus\fa_0\rightarrow \fa\oplus\fa_0$ is antisymmetric and satisfies $D_\fa\circ\hat U=\hat U\circ D_\fa$ and $\theta_\fa\circ\hat U=- \hat U\circ\theta_\fa$ if and only if $\hat U(A_0)=0$ and $\hat U|_{\fa_0}\in\fso(\fa_0)^D_-$. As in (\ref{ESym}) we have  $\fso(\fa_0)^D_-\cong \Sym^2((\fa_0)_-^*).$ Let $b\in\Sym^2((\fa_0)_-^*)$ be the element corresponding to $\hat U$ under this isomorphism. The cochain $\hat\tau\in C^1(\fl,\fa)$ satisfies Condition $(ii)$ in Prop.~\ref{Pphi} if and only if $\hat \tau(X)=\hat a$ and $\hat \tau(Y)=D_\fa(\hat a)$ for some $\hat a \in (\fa_0)_-$. The remaining conditions in Prop.~\ref{Pphi} are trivially satisfied. This gives us an isomorphism
\begin{eqnarray}\label{ER2}
 \Out(\fg)^D_-&\longrightarrow &  \RR\oplus \Sym^2((\fa_0)_-^*) \oplus (\fa_0)_-\nonumber\\
 \ph(\hat S,\hat U,\hat \tau,0)&\longmapsto & (\mu,b,\hat a).
 \end{eqnarray}
Tensoring by $R$ gives the assertion of the proposition.
 \qed

Abbreviating our notation we will write $\Aut(\fg)$ for the automorphism group of the extrinsic symmetric triple $\fg$.  In the following proposition we use the notation from Definition~\ref{DF} for $\fg=\fd_{\alpha,\gamma}(\fl,\Phi_\fl,\fa\oplus\fa_0)$.
\begin{pr}
The automorphism group $\Aut(\fg)$ equals
$$\left\{F(S,U,\tau,\sigma)\left| \begin{array}{l}U(A_0)=A_0,\ U((\fa_0)_-)=(\fa_0)_-,\ U|_{(\fa_0)_+}=-D_\fa U|_{(\fa_0)_-} D_\fa,
\\[0.2ex] S=\pm \Id,\ \tau(Y)=:\bar a\in (\fa_0)_-,\ \tau(X)=-D_\fa(\bar a),\ \sigma=0\end{array}\right.\right\}.$$
In particular, it
is isomorphic to the semidirect product $(\ZZ_2\times \grO((\fa_0)_-))\ltimes (\fa_0)_-$. 
\end{pr} 
\proof Let $F(S,U,\tau,\sigma)$ be an automorphism of $\fd_{\alpha,\gamma}(\fl,\Phi_\fl,\fa\oplus\fa_0)$. Since $\theta_\fl S=S\theta_\fl$, the isomorphism $S:\fl\rightarrow \fl$ must satisfy $S(X)=\lambda X$, $S(Y)=\mu Y$. The condition $D_\fl S=S D_\fl$ now implies $\lambda=\mu$. Because of $D_\fa U=UD_\fa$ we have $U(A_0)=\kappa A_0$, $\kappa=\pm1$. Thus also $U(\fa_0)=\fa_0$, which implies $U((\fa_0)_-)=(\fa_0)_-$ and $U|_{(\fa_0)_+}=-D_\fa U|_{(\fa_0)_-} D_\fa$ by $D_\fa$-invariance of $U$.
Moreover, $d\tau=0$ implies
$$A_0=\alpha(X,Y)=(U\circ S^*\alpha)(X,Y)=\kappa\lambda^2 A_0,$$
thus $1=\kappa\lambda^2=\kappa$. Finally, $\sigma(X,Y)=(\theta_\fl^*\sigma)(X,Y)=-\sigma(X,Y)$ implies $\sigma=0$. Conversely, if $S,U,\tau,\sigma$ are given as in the proposition, then one easily checks that $F(S,U,\tau,\sigma)$ is an automorphism.

Since $U$ is defined by  $\bar U:=U|_{(\fa_0)_-}$ and $\tau$ is defined by $\tau(Y)=:\bar a\in(\fa_0)_-$ the map 
\begin{eqnarray}\label{EFs}
\Aut(\fg)&\longrightarrow & (\ZZ_2\times \grO((\fa_0)_-))\ltimes (\fa_0)_-\nonumber\\
F(\lambda\cdot\Id,U^{-1},\lambda U^{-1}\circ\tau,0)&\longmapsto&(\lambda,\bar U,\bar a)
\end{eqnarray}
is bijective. Moreover, denoting the preimage $F(\lambda\cdot\Id,U^{-1},\lambda U^{-1}\circ\tau,0)$ of $(\lambda,\bar U,\bar a)$ by $f(\lambda,\bar U,\bar a)$ we obtain $$f(\lambda_1,\bar U_1,\bar a_1)\cdot f(\lambda_2,\bar U_2,\bar a_2)=f(\lambda_1\lambda_2,\bar U_1 \bar U_2,\bar a_1 +\lambda_1 \bar U_1\bar a_2).$$

\vspace{-4ex}\qed

We consider $(\fa_0)_-$ as a Euclidean space and denote by $\Iso((\fa_0)_-):=\grO((\fa_0)_-)\ltimes (\fa_0)_-$ its affine isometry group. Elements of this group are denoted by $(\bar U,\bar a)$, where $\bar U\in \grO((\fa_0)_-)$ and $\bar a\in(\fa_0)_-$. As usual $(\cdot)^\sharp:\fa_0^*\otimes R\rightarrow\fa_0\otimes R$ will denote the isomorphism defined by $\ip_\fa$.
\begin{pr} \label{PH2Aut} There is a bijection from
$H^2(\fg,R)_-^D/(\Aut(\fg)\times \GL(R))$ to
$$\left(R\oplus S^2((\fa_0)_-,R)\oplus (\fa_0)_-\otimes R\right)/(\Iso((\fa_0)_-)\times\GL(R)),$$
where $(\bar U,\bar a)\in \Iso((\fa_0)_-)$ acts by
$$(r_0,B,\eta)*(\bar U,\bar a)=\left( r_0,\,\bar U^*B,\, \bar U^{-1}( \eta- B(\bar a,\cdot\,)^\sharp-\bar a\otimes r_0)\right),$$
and $g\in\GL(R)$ acts on each summand in a natural way. 
\end{pr}
\proof
First we determine the action of $\Aut(\fg)$ on $H^2(\fg,\RR)^D_-\cong \Out(\fg)^D_-$. Let $\phi(\mu,b,\hat a)$ be the preimage of $(\mu,b,\hat a)$ under the map (\ref{ER2}). Then 
$$f(\lambda,\bar U,\bar a)^{-1}\cdot \phi(\mu,b,\hat a)\cdot f(\lambda,\bar U,\bar a)\equiv \phi(\mu,\bar U^* b,\lambda\bar U^{-1}(\hat a- b(\bar a,\cdot\,)^\sharp -\mu \bar a))\ \mod \ad(\fg^+_-).$$
Now tensoring by $R$ gives the asserted formula for the action.

In particular, the element $(-1,-\Id,0)\in(\ZZ_2\times \grO((\fa_0)_-))\ltimes (\fa_0)_-$ acts as identity. Hence the orbits of $(\ZZ_2\times \grO((\fa_0)_-))\ltimes (\fa_0)_-$ are equal to the orbits of $\grO((\fa_0)_-)\ltimes (\fa_0)_-=\Iso((\fa_0)_-)$.
\qed
\begin{de}\label{Dind}
An element $(r_0,B,\eta)\in R\oplus S^2((\fa_0)_-,R)\oplus (\fa_0)_-\otimes R$ is called decomposable if and only if there are decompositions $R= R'\oplus R''$ and $(\fa_0)_-=\fa'\oplus\fa''$, $\fa'\perp\fa''$ with $\fa''\oplus R''\not=0$ such that  
\begin{itemize}
\item[(i)] $r_0\in R'$;
\item[(ii)] $B(\fa',\fa')\subset  R'$, $B(\fa'',\fa'')\subset R''$, $B(\fa',\fa'')=0$;
\item[(iii)] $\exists\,\eta'\in \fa'\otimes R' \ \exists\, a''\in\fa'': \eta=\eta'+ a''\otimes r_0+ B(a'',\cdot\,)^\sharp.$
\end{itemize}
Otherwise $(r_0,B,\eta)$ is called indecomposable. 
\end{de}
Decomposability of $(r_0,B,\eta)$ only depends on the $\Aut(\fa_0)\times \GL(R)$-orbit $[r_0,B,\eta]$ of~$(r_0,B,\eta)$, hence we can speak of decomposable and indecomposable orbits.
\begin{pr}
\quad The weak extrinsic symmetric triple corresponding to the orbit $[r_0,B,\eta]$ is full and indecomposable if and only if $[r_0,B,\eta]$ is indecomposable.
\end{pr}
\proof We apply Prop.~\ref{Plt}. In our situation the first and the third condition of Prop.~\ref{Plt} are obviously satisfied. Let us turn to the second condition. We reformulate decomposability for $[\omega]\in H^2(\fg,R)^D_-$ in terms of the element $(r_0,B,\eta)$  corresponding to $[\omega]$ under the isomorphism $H^2(\fg,R)_-^D\cong R\oplus S^2 ((\fa_0)_-,R)\oplus (\fa_0)_-\otimes R$. The cohomology class $[\omega]$ is in the $\Aut(\fg)$-orbit of a class $[\omega']+[\omega'']$ as described in the second conditon of Prop.~\ref{Plt} if and only if
\begin{eqnarray}
(r_0,B,\eta)&=&(r,B'+B'',\eta')*(\bar U,\bar a)\nonumber\\
&=&\left(r,\,\bar U^*(B'+B''),\, \bar U^{-1}(\eta'-(B'+B'')(\bar a,\cdot\,)^\sharp-\bar a\otimes r )\right),
\label{Edecomp}
\end{eqnarray}
where $r\in R'$, $\eta'\subset\fa'\otimes R'$ and
$$
B'(\fa',\fa')\subset R', \quad B'(\fa'',(\fa_0)_-)=0,\quad
B''(\fa'',\fa'')\subset R'', \quad B''(\fa',(\fa_0)_-)=0
$$
for $\fa':=(\fa_0')_-$ and $\fa'':=(\fa_0'')_-$. Assume that $\fg$ is decomposable. Then Conditions $(i)$ -- $(iii)$ of Def.~\ref{Dind} hold for the decompositions $R=R'\oplus R''$ and $(\fa_0)_-=\bar U^{-1}(\fa')\oplus \bar U^{-1}(\fa'')$. Indeed, by (\ref{Edecomp})
\begin{eqnarray*}
&&r_0= r\ \in\ R'\\
&&B( \bar U^{-1}(\fa'),  \bar U^{-1}(\fa'))\ =\ (B'+B'')(\fa',\fa')\ \subset\ R'\\
&&B( \bar U^{-1}(\fa''),  \bar U^{-1}(\fa''))\ =\ (B'+B'')(\fa'',\fa'')\ \subset\ R''\\
&&B( \bar U^{-1}(\fa'),  \bar U^{-1}(\fa'')) \ =\ (B'+B'')(\fa',\fa'')\ =\ 0.
\end{eqnarray*}
Furthermore, if we decompose $\bar a=\bar a'+\bar a''$ with respect to $\fa'\oplus\fa''$, then 
\begin{eqnarray*} 
\bar U^{-1}(\eta'-(B'+B'')(\bar a,\cdot\,)^\sharp-\bar a\otimes r)&=&\bar U^{-1}(\eta'-B'(\bar a,\cdot\,)^\sharp-\bar a'\otimes r)-\bar U^{-1}(\bar a'')\otimes r\\
&&-B(\bar U^{-1}(\bar a''),\cdot\,)^\sharp
\end{eqnarray*}
 and $\bar U^{-1}(\eta'-B'(\bar a,\cdot\,)^\sharp-\bar a'\otimes r)\subset \bar U^{-1}(\fa')\otimes R'$.

 Conversely, let $(r_0,B,\eta)$ satisfy Conditions $(i)$ -- $(iii)$ of Def.~\ref{Dind} for decompositions $R=R'\oplus R''$ and $\fa=\fa'\oplus\fa''$. Define  $B':= B\circ\proj_{\fa'}$, $B'':=B\circ\proj_{\fa''}$. Then
\begin{eqnarray*}
(r_0,B,\eta)&=& (r_0,B'+B'',\eta'+a''\otimes r_0+(B'+B'')(a'',\cdot)^\sharp) \\
&=& (r_0,B'+B'',\eta') *(\Id,-a''),
\end{eqnarray*}
hence $(r_0,B,\eta)$ satisfies (\ref{Edecomp}). Thus $\fg$ is decomposable.\qed
\paragraph{The case \fbox{$\fl=\fh(1)$}}
Let $\fg_1:=\dd$ be as in Thm.~\ref{Tfull}, \ref{I3}.~defined.   

\begin{lm} For $\fl=\fh(1)$ and $\fg=\fg_1\oplus\fa_0$ with $\fg_1:=\dd$ as considered above 
$$H^2(\fg,R)_-^D\,\cong\, R\,\oplus \,S^2((\fa_0)_-,R)\,\oplus\, (\fa_0)_-\otimes R.$$
\end{lm}
\proof First let us check which $\ph(\hat S,\hat U,\hat \tau,\hat \sigma)$ are in $\Der(\fg)^D_-$. The linear map $\hat S:\fl\rightarrow\fl$ is a derivation satisfying $D_\fl\circ\hat S=\hat S\circ D_\fl$ and $\theta_\fl\circ\hat S=- \hat S\circ \theta_\fl$ if and only if $\hat S(X)=\mu Y$, $\hat S(Y) =-\mu X$ and $S(Z)=0$ for some $\mu\in\RR$. The cochain $\hat\tau\in C^1(\fl,\fa)$ satisfies Condition (ii) in Prop.~\ref{Pphi} if and only if $\hat \tau(X)=\hat a$, $\hat \tau(Y)=D_\fa(\hat a)$ and $\hat \tau(Z)=0$ for some $\hat a \in \RR\cdot A_1\oplus (\fa_0)_-$. 
Furthermore, $\hat U:\fa\oplus\fa_0\rightarrow \fa\oplus\fa_0$ is antisymmetric and satisfies $D_\fa\circ\hat U=\hat U\circ D_\fa$ and $\theta_\fa\circ\hat U=- \hat U\circ\theta_\fa$ as well as 
$$\hat U(\alpha(L_1,L_2))+\alpha(\hat S L_1,L_2)+\alpha(L_1,\hat S L_2)=-d\hat \tau (L_1,L_2)=0$$ 
for all $L_1,L_2\in\fl$ if and only if 
$\hat U(A_1)=-\mu A_2,\ \hat U(A_2)=\mu A_1$ and $\hat U|_{\fa_0}\in\fso(\fa_0)^D_-$. As in (\ref{ESym}) we have  $\fso(\fa_0)^D_-\cong \Sym^2((\fa_0)_-^*).$ Let $b\in\Sym^2((\fa_0)_-^*)$ be the element corresponding to $\hat U$ under this isomorphism. Now we consider $\hat \sigma:\fl\rightarrow \fl^*$.  Then $\hat \sigma$ satisfies Condition (ii) in Prop.~\ref{Pphi} if and only if $\hat\sigma(X)=s\sigma_Y$, $\hat\sigma(Y)=-s\sigma_X$ and $\hat\sigma(Z)=0$ for some $s\in\RR$.
The remaining conditions in Prop.~\ref{Pphi} are trivially satisfied. 

Now we determine $\ad(\fg^+_-)$. Here we have $\fg^+_-=\Span\{Z,\sigma_Z\}$ and
$$[Z,X]=-A_1,\ [Z,Y]=-A_2,\ [\sigma_Z,X]=\sigma_Y,\ [\sigma_Z,Y]=-\sigma_X.$$
 This gives an isomorphism
\begin{eqnarray*}
\Out(\fg)^D_-&\longrightarrow & \RR\oplus \Sym^2((\fa_0)_-^*) \oplus (\fa_0)_-\nonumber\\
 \ph(\hat S,\hat U,\hat \tau,\hat \sigma)+\ad(\fg^+_-)& \longmapsto & (\mu,b,\proj_{\fa_0}\hat a).
\end{eqnarray*}
Tensoring by $R$ gives the assertion of the proposition.
 \qed

\begin{pr}
The automorphism group $\Aut(\fg)$ of the  extrinsic symmetric triple $\fg=\fg_1\oplus\fa_0$ equals
$$\left\{F(S,U,\tau,\sigma)\left| \begin{array}{l}S(X)=\lambda X, \ S(Y)=\lambda Y,\, S(Z)=Z,\\ U(A_1)=\lambda A_1,\  U(A_2)=\lambda A_2,\ \lambda=\pm1,\\ U((\fa_0)_-)=(\fa_0)_-,\ U|_{(\fa_0)_+}=-D_\fa U|_{(\fa_0)_-} D_\fa,
\\ \tau(Y)=:\bar a\in (\fa_0)_-,\ \tau(X)=-D_\fa(\bar a),\ \tau(Z)=0,\ \sigma=0\end{array}\right.\right\}.$$
In particular, it
is isomorphic to the semidirect product $(\ZZ_2\times \grO((\fa_0)_-))\ltimes (\fa_0)_-$. 
\end{pr} 
\proof Let $F(S,U,\tau,\sigma)$ be an automorphism of $\fd_{\alpha,\gamma}(\fl,\Phi_\fl,\fa\oplus\fa_0)$. Since $\theta_\fl S=S\theta_\fl$ and $D_\fl S=SD_\fl$, the isomorphism $S:\fl\rightarrow \fl$ must satisfy $S(X)=\lambda X$, $S(Y)=\lambda Y$ and $S(Z)=\mu(Z)$ for some $\mu,\lambda\in\RR$. Because of $[SX,SY]=S[X,Y]$ we get $\mu=\lambda^2$. Condition $(iii)$ of Prop.~\ref{PF} gives 
 $$\lambda^3 U(A_1)=(U\circ S^*\alpha)(X,Z)=\alpha(X,Z)-d\tau(X,Z)=A_1.$$
 Since $U$ is an isometry we get $\lambda=\pm1$. Moreover, $D_\fa U=UD_\fa$ implies $U(A_2)=\lambda A_2$. Now we see $U((\fa_0)_-)=(\fa_0)_-$ and $U|_{(\fa_0)_+}=-D_\fa U|_{(\fa_0)_-} D_\fa$. Finally,  $$0=(U\circ S^*\alpha)(X,Y)=\alpha(X,Y)-d\tau(X,Y)=-d\tau(X,Y)=\tau([X,Y])=\tau(Z)$$
 and $\sigma=0$ because of $\sigma(X,Y)=(\theta_\fl^*\sigma)(X,Y)=-\sigma(X,Y)$ and $0=\sigma(D_\fl L, Z)+\sigma(L, D_\fl Z)=\sigma(D_\fl L, Z)$ for $L\in\{X,Y\}$.
 
 For $\lambda=\pm1$ let the isomorphism $S_\lambda:\fl\rightarrow\fl$ be defined by $S_\lambda(X)=\lambda X$, $S_\lambda(Y)=\lambda Y$, $S_\lambda(Z)=Z$. As in the case $\fl=\RR^2$ we identify $\Aut(\fg)$ with $(\ZZ_2\times \grO((\fa_0)_-))\ltimes (\fa_0)_-$.  Here we use the bijection 
 \begin{eqnarray}\label{EFs2}
\Aut(\fg)&\longrightarrow & (\ZZ_2\times \grO((\fa_0)_-))\ltimes (\fa_0)_-\nonumber\\
F(S_\lambda,U^{-1},\lambda U^{-1}\circ\tau,0)&\longmapsto&(\lambda,\bar U,\tau(Y)),
\end{eqnarray}
where $\bar U:= U|_{(\fa_0)_-}$.
 \qed
 
 In the same way as for $\fl=\RR^2$ we obtain:
 \begin{pr} There is a bijection between 
$H^2(\fg,R)_-^D/(\Aut(\fg)\times \GL(R))$ and
$$\left(R\oplus S^2((\fa_0)_-,R)\oplus  (\fa_0)_-\otimes R\right)/(\Iso((\fa_0)_-)\times\GL(R)),$$
where $\Iso((\fa_0)_-)\times\GL(R)$ acts on $R\oplus S^2((\fa_0)_-,R)\oplus  (\fa_0)_-\otimes R$ as defined in Prop.~\ref{PH2Aut}.
\end{pr}
 
In particular, we have the same condition for decomposabiltity as for $\fl=\RR^2$. 
\begin{pr}
\quad The weak extrinsic symmetric triple corresponding to the orbit $(r_0,\eta,B)\cdot (\Iso((\fa_0)_-)\times\GL(R))$ is indecomposable if and only if $(r_0,\eta,B)$ is indecomposable in the sense of Definition~\ref{Dind}.
\end{pr} 
 \paragraph{The case \fbox{$\fl\in\{\fsl(2,\RR),\, \fsu(2)\}$}}

Let $\fg_1:=\dd$ be as in Thm.~\ref{Tfull}, \ref{I4}.\,or \ref{I5}.\,defined. We introduce the notation
$$\frak V:=({\fa_3})^k,\ \hat \frak V:=(\hat {\fa}_3)^l,\ \frak W:=({\fa_4})^m.$$
Moreover, $V:=\frak V\cap \fa_-^-$, $\hat V:=\hat\frak V\cap \fa_-^-$, $W:=\frak W\cap \fa_-^-$.
\begin{lm} For $\fl\in\{\fsl(2,\RR),\, \fsu(2)\}$ and $\fg=\fg_1\oplus\fa_0$ with $\fg_1:=\dd$ as considered above we have
$$
H^2(\fg,R)_-^{D}\cong  \Hom(V\otimes\hat V,R)\oplus S^2(W,R)\oplus S^2((\fa_0)_-,R). 
$$
\end{lm} 
\proof Note first that in this case Condition $(iii)$ in Prop.~\ref{Pphi} reduces to the equation $d\hat\tau=0$, which is equivalent to $\hat \tau =dA$ for some $A\in\fa\oplus\fa_0$ since $\fl$ is semisimple. Therefore we will assume that $\hat\tau=dA$ holds and then confine ourselves to conditions $(i)$ and $(ii)$.

Now we show that we can identify 
$$\Out(\fg)^D_-\cong \{\ph(\hat S,\hat U,\hat \tau,\hat \sigma)\in \Der(\fg)^D_-\mid \hat S=0,\ \hat \tau=0,\ \hat \sigma=0\}.$$
Any derivation of $\fl$ that satisfies the conditions in $(i)$ is a multiple of $\ad(X)$. Hence for any such $\hat S$ there is an inner derivation $\ph\in\ad(\fg^+_-)$ such that $\ph=\ph(\hat S,\hat U,0,\hat \sigma)$ for some $\hat U$ and $\hat \sigma$. The map $\hat \tau=dA$ satisties the conditions in $(ii)$ if and only if $A\in\fa_-^+$. For $\hat \tau =dA$ with $A\in\fa_-^+$ there is a derivation $\ph\in \ad(\fg^+_-)$ such that $\ph=\ph(0,0,\hat \tau,0)$, namely $\ph=-\ad(A)$. The only linear maps $\hat\sigma:\fl\rightarrow \fl^*$ that satisfy the conditions in $(ii)$ are multiples of $\ad (\sigma_X)|_{\fl}$. Thus for all $\hat\sigma$ satisfying $(ii)$ there exists a  derivation $\ph\in \ad(\fg^+_-)$ with $\ph=\ph(0,0,0,\hat\sigma)$.  We conclude, that any outer derivation can be represented by a map $\ph(0,\hat U,0,0)$. Moreover, if $\ph=\ph(0,\hat U,0,0)$ is in $\ad(\fg_+^-)$, then $\ph=0$.

It remains to check which $\ph(0,\hat U,0,0)$ are in $\Der(\fg)^D_-$. This is equivalent to check which $\hat U\in \fso(\fa\oplus\fa_0)$ commute with $D_\fa$, anticommute with $\theta_\fa$ and satisfy $[\rho(L),\hat U]=0$ for all $L\in\fl$. The latter equation implies $\hat U(\frak V\oplus \hat\frak V)\subset \frak V\oplus \hat\frak V$, $\hat U(\frak W)\subset \frak W$ and $\hat U(\fa_0)\subset\fa_0$. Moreover, it yields
$$U_1:=\hat U|_V:V=\frak V \cap \Ker H \longrightarrow (\frak V\cap \hat\frak V) \cap \Ker H.$$ Since $D_\fa\circ\hat U=\hat U\circ D_\fa$ and $\theta_\fa\circ\hat U=- \hat U\circ\theta_\fa$ we get on the other hand 
$$U_1:V=\frak V \cap \fa_-^- \longrightarrow  (\frak V\oplus \hat\frak V) \cap \fa_+^-,$$
thus $U_1(V)\subset  \hat\frak V \cap \fa_+^-$ because of $\frak V\cap \Ker H\cap \fa_+^-=0$. Conversely, if there is given a linear map $U_1:V \rightarrow \hat\frak V \cap \fa_+^-$, then we can extend $U_1$ in a unique way to a linear map $\hat U:\frak V\rightarrow \hat\frak V$ such that $[\rho(L)|_{\frak V},\hat U]=0$ holds for all $L\in\fl$. Moreover, we can $\hat U$ further extend to an antisymmetric map $\hat U\in\fso( \frak V\oplus \hat\frak V)$. This extension is also unique. It automatically satisfies $D_\fa\circ\hat U=\hat U\circ D_\fa$ and $\theta_\fa\circ\hat U=- \hat U\circ\theta_\fa$. In particular, $\hat U|_{\frak V\oplus\hat \frak V}$ can be identified with the bilinear map $b_1:= \langle U_1(\cdot),D(\cdot) \rangle \in V^*\otimes \hat V^*$.

In a similar way we treat $\hat U|_{\frak W}$. This map is completely determined by $U_2 :=\hat U|_W:W\rightarrow \frak W\cap \fa_+^-$. The antisymmetry of $\hat U$ implies 
\begin{equation}\label{EU2}
\langle U_2w_1,D w_2\rangle =\langle U_2 w_2, D w_1\rangle
\end{equation} 
for all $w_1,w_2\in W$. Conversely, any linear map $U_2 :W\rightarrow \frak W\cap \fa_+^-$ that satisfies (\ref{EU2}) defines uniquely an antisymmetric map $\hat U:\frak W\rightarrow \frak W$ for which $[\rho(L)|_{\frak W},\hat U]=0$ holds for all $L\in\fl$. This map automatically satisfies $D_\fa\circ\hat U=\hat U\circ D_\fa$ and $\theta_\fa\circ\hat U=- \hat U\circ\theta_\fa$. Hence $\hat U|_{\frak W}$ can be identified with the symmetric  bilinear map $b_2= \langle U_2 (\cdot),D(\cdot) \rangle\in \Sym^2 (W^*)$.

Finally, as we already know, $\hat U|_{\fa_0}\in\fso(\fa_0)$ corresponds to an element $b\in \Sym^2((\fa_0)_-^*)$. 

This gives us an isomorphism
\begin{eqnarray*} 
\Out(\fg)^D_-&\longrightarrow&  (V^*\otimes \hat V^*) \otimes \Sym^2 (W^*)\otimes \Sym^2((\fa_0)_-^*)\\
\ph=\ph(0,\hat U,0,0)&\longmapsto& \left(\ b_1,\ b_2,\ b\ \right).
\end{eqnarray*}
Tensoring by $R$ yields the assertion.
\qed  
\begin{pr}
The automorphism group of the  extrinsic symmetric triple $\fg\oplus\fa_0$ consists of all $F(S,U,\tau,\sigma)$ satisfying
$$ \begin{array}{l}S(X)= X, \ S(Y)=\lambda Y,\, S(H)=\lambda H,\ \lambda=\pm1,\\ U((\fa_0)_-)=(\fa_0)_-,\ U|_{(\fa_0)_+}=-D_\fa U|_{(\fa_0)_-} D_\fa,\\
U(V)=V,\ U(\hat V)=\hat V,\ U(W)=W,\\
U|_{\frak a^-_+}=D_\fa U|_{\fa^-_- }D_\fa^{-1} ,\  U|_{\fa^+_+}=\lambda Y\circ U|_{\fa^-_-} \circ Y^{-1},\ U|_{\fa^+_-}=\lambda H\circ U|_{\fa_-^-} \circ H^{-1} , \\ \tau=dA \mbox{ for some fixed } A\in\fa^+_+,\ \mbox{and } \sigma=0.\end{array}$$
In particular, it
is isomorphic to $\Big((\ZZ_2\times \grO(V)\times\grO(\hat V)\times\grO(W))\ltimes \fa_+^+\Big)\times \grO((\fa_0)_-)$. 
\end{pr} 
\proof First note that $S(X)=\mu(X)$, $S(Y)=\lambda Y$, $S(H)=\lambda(H)$ holds for some $\mu, \lambda\in\RR$ since $S$ commutes with $D_\fl$ and $\theta_\fl$. Since $S$ is a Lie algebra isomorphism $\mu=\lambda^2$ must hold.  We define the linear map $S_\lambda:\fl\rightarrow \fl$ by $S_\lambda(X)=X$, $S_\lambda(Y)=\lambda Y$ and  $S_\lambda(H)=\lambda H$. If $U:\fa\oplus\fa_0\rightarrow \fa\oplus\fa_0$ is an isometry that satisfies $\rho(S_\lambda(\cdot))\circ U=U\circ \rho(\cdot)$, then $U$ must respect the decomposition $\fa=\frak V\oplus\hat \frak V\oplus\frak W\oplus\fa_0$. Moreover, if we take into account that $U$ commutes with  $D$ and $\theta$, we see that $U$ must have the properties claimed in the proposition.  Conversely, if $U$ satisfies the conditions in the proposition, then $(S_\lambda,U^{-1}): ((\fl,\Phi_\fl),\fa\oplus\fa_0)\rightarrow ((\fl,\Phi_\fl),\fa\oplus\fa_0)$ is an isomorphism of pairs. 

Since $\alpha=0$ we have $d\tau=0$. Because of $H^1(\fl,\fa\oplus\fa_0)=0$ this implies $\tau =d A$ for some $A\in\fa$. Now $\tau\in C^1(\fl,\fa)^{(D,\theta)}$ yields $A\in\fa^+_+$.  Moreover, $\sigma=0$ since  $C^2(\fl)^{(D,\theta)}=0$. Obviously, $(S_\lambda,U)^*(0,\gamma)=(0,\gamma)(\tau,0)^{-1}$ holds for any $S_\lambda, U,\tau$ as chosen above. 

For any automorphism of pairs $(S_\lambda,U^{-1}): ((\fl,\Phi_\fl),\fa\oplus\fa_0)\rightarrow ((\fl,\Phi_\fl),\fa\oplus\fa_0)$ we define 
$$\bar U_1:= U|_{\fa^-_-}\in\grO(V)\times\grO(\hat V)\times\grO(W)\subset \grO(\fa_-^-),\quad \bar U_2:= U|_{(\fa_0)_-}\in  \grO((\fa_0)_-).$$
Then
\begin{eqnarray}\label{fast}
\Aut(\fg\oplus\fa_0)&\longrightarrow & \Big((\ZZ_2\times \grO(V)\times\grO(\hat V)\times\grO(W))\ltimes \fa_+^+\Big)\times \grO((\fa_0)_-)\nonumber \\ F(U^{-1},U^{-1}dA\circ S_\lambda, S_\lambda,0)&\longmapsto& (\lambda,\bar U_1,A,\bar U_2)
\end{eqnarray}
is a group isomorphism.
Here  $(\lambda,\bar U_1)\in\ZZ_2\times \grO(V)\times\grO(\hat V)\times\grO(W)$ acts on $\fa^+_+$ by $\lambda Y\circ \bar U_1\circ Y^{-1}$.
\qed
\begin{pr} There is a bijection from
$H^2(\fg,R)_-^D/(\Aut(\fg)\times \GL(R))$ to
$$\Big({}^{\mbox{\begin{small}$\Hom(V\otimes\hat V,R)$\end{small}}}/_{\mbox{\begin{small}$\grO(V)\times \grO(\hat V)$\end{small}}} \oplus {}^{\mbox{\begin{small}$S^2(W,R)$\end{small}}}/{}_{\mbox{\begin{small}$\grO(W)$\end{small}}} \oplus {}^{\mbox{\begin{small}$S^2((\fa_0)_-,R)$\end{small}}}/_{\mbox{\begin{small}$\grO((\fa_0)_-)$\end{small}}} \Big)/\mbox{\begin{small}$\GL(R)$\end{small}},$$
where $\GL(R)$ and all orthogonal groups act in the natural way.
\end{pr}
\proof Take $(b_1,b_2,b)\in (V^*\otimes \hat V^*) \otimes \Sym^2 (W^*)\otimes \Sym^2((\fa_0)_-^*)$ and consider the corresponding $\ph=\ph(0,\hat U,0,0) \in\Out(\fg^+_-)$. Then the automorphism  $F:=F(U^{-1},U^{-1}dA\circ S_\lambda, S_\lambda,0)$ that corresponds to $(\lambda,\bar U_1,A,\bar U_2)$ according to (\ref{fast}) acts by conjugation on $\ph$. One easily computes $F^{-1}\ph F\equiv \ph(0,U^{-1}\hat U U,0,0) \mod \ad(\fg^+_-)$, which translates into  
$$
(b_1,b_2,b)*(\lambda,\bar U_1,A,\bar U_2) = (\bar U_1^* b_1, \bar U_1^* b_2, \bar U_2^*b)
$$
 Taking the tensor product by $R$ gives the assertion. \qed

\begin{de} An element $(B_1,B_2,B)\in \Hom(V\otimes\hat V,R)\times S^2(W,R)\times S^2((\fa_0)_-,R)$ is called decomposable if there 
  are decompositions $R= R'\oplus R''$ and $(\fa_0)_-=\fa'\oplus\fa''$, $\fa'\perp\fa''$ with $\fa''\oplus R''\not=0$ such that  
\begin{itemize}
\item[(i)] $B_1(V,\hat V)\subset R'$, $B_2(W,W)\subset R'$; and
\item[(ii)] $B(\fa',\fa')\subset  R'$, $B(\fa'',\fa'')\subset R''$, $B(\fa',\fa'')=0$.
\end{itemize}
\end{de}
Decomposability of $B$ only depends on the $\Aut(\fa_0)\times \GL(R)$-orbit of~$B$, hence we can speak of decomposable and indecomposable orbits.
\begin{pr}
The weak extrinsic symmetric associated to the orbit of $(B_1,B_2,B)$ is full and indecomposable if and only if $(B_1,B_2,B)$ is indecomposable.
\end{pr}

Summarizing the results of this section we obtain a classification of extrinsic symmetric spaces of Lorentz type. Recall that $\fa_0$ denotes an abelian extrinsic symmetric triple $(\fa_0,\Phi_\fa,\ip_\fa)$, $\Phi_\fa=(D_\fa,\theta_\fa)$, for which $\ip_\fa$ is positive definite and  $D_\fa^2=-\Id$.
\begin{theo}\label{Tende}
Each full weak extrinsic symmetric triple of Lorentz type that does not have simple ideals is isomorphic to a normal extension of an extrinsic symmetric triple $\dd\oplus\fa_0$ by a vector space $R$, where $\dd$ is one of the full Lorentzian extrinsic symmetric triples from Thm.~\ref{Tfull}.  More exactly:
\begin{itemize}
\item[(i)] Full and indecomposable weak extrinsic symmetric triples of Lorentz type that are normal extensions of $\tilde\fa:=\RR^{2,0}\oplus\fa_0$  by $R$ are in bijection to indecomposable elements of the orbit space 
$$ S^2(\tilde\fa_-,R)/(\grO(\tilde\fa_-)\times\GL(R)).$$
\item[(ii)] Full and indecomposable weak extrinsic symmetric triples of Lorentz type that are normal extensions by $R$ of $\dd\oplus\fa_0$ for $\dd$ as in Thm.~\ref{Tfull}, case \ref{I2}  or case~\ref{I3}, are in bijection to indecomposable elements of the orbit space 
$$ \left(R\oplus S^2((\fa_0)_-,R)\oplus  (\fa_0)_-\otimes R\right)/(\Iso((\fa_0)_-)\times\GL(R)).$$
\item[(iii)] Full and indecomposable weak extrinsic symmetric triples of Lorentz type that are normal extensions by $R$ of $\dd\oplus\fa_0$ for $\dd$ as in Thm.~\ref{Tfull}, case \ref{I4} or case~\ref{I5}, are in bijection to  indecomposable elements of the orbit space 
$$\Big({}^{\mbox{\begin{small}$\Hom(V\otimes\hat V,R)$\end{small}}}/_{\mbox{\begin{small}$\grO(V)\times \grO(\hat V)$\end{small}}} \oplus {}^{\mbox{\begin{small}$S^2(W,R)$\end{small}}}/{}_{\mbox{\begin{small}$\grO(W)$\end{small}}} \oplus {}^{\mbox{\begin{small}$S^2((\fa_0)_-,R)$\end{small}}}/_{\mbox{\begin{small}$\grO((\fa_0)_-)$\end{small}}} \Big)/\mbox{\begin{small}$\GL(R)$\end{small}}.$$
Here $V=({\fa_3})^k\cap \fa_-^-,\ \hat  V=(\hat {\fa}_3)^l\cap \fa_-^-,\  W=({\fa_4})^m\cap \fa_-^-.$
\end{itemize}
\end{theo}
\begin{ex}{\rm We want to compare our classification of extrinsic symmetric spaces $M\hookrightarrow V$ in the case $\dim M=2$, $V=\RR^{1,3}$ with the classification of parallel surfaces in $\RR^{1,3}$ from \cite{CV}. Any non-degenerate surface  is full either in a non-degenerate subspace or in the orthogonal sum $\RR^2\oplus R$ of $\RR^2$ with the standard scalar product and a one-dimensional isotropic space $R$, where the latter case can only occur in the case of a spacelike surface. 

First let $M$ be spacelike. If $M$ is full in a non-degenerate subspace, then Prop.~\ref{Pqr} shows that the associated extrinsic symmetric triple must be reductive. All such $M$ are listed in \cite{CV}, Thm.\,7.1., (1) -- (6). If $M$ is full in $\RR^2\oplus R$, then the associated weak extrinsic symmetric triple is a normal extension of the four-dimensional abelian extrinsic symmetric triple $\fa$ associated with the embedding $\RR^2\hookrightarrow \RR^2$ by the one-dimensional space $R$. By  Prop.~\ref{PRt} the possible isometry classes of $M\hookrightarrow\RR^2\oplus R$ correspond to the elements of
$$S^2(\fa_-,R)/(\grO(\fa_-)\times\GL(R))=\left\{\left.\left(\begin{array}{cc}1&0\\0&\lambda\end{array}\right) \right| \lambda\in\RR,\ \lambda=0 \mbox{ or } |\lambda|>1   \right\}.$$
The associated surfaces have the parametrisation $(u,v)\mapsto (u,v,u^2+\lambda v^2)\in\RR^2\oplus R$. For $\lambda=-1$ this is case (7) in \cite{CV}, $\lambda=1$ corresponds to case (8) and for the remaining $\lambda$ we are in case (9). Note that in \cite{CV}, case (9) one should replace $b\in\RR$ by $b>0$ to ensure uniqueness.

Now let $M$ be timelike. Then $M$ is full in some non-degenerate subspace of $\RR^{1,3}$. If the associated extrinsic symmetric triple is reductive then $M$ must be one of the spaces in \cite{CV}, (1) -- (5). Otherwise it is the space in Corollary \ref{Co}, (\ref{II2}), which corresponds to  \cite{CV}, (6).
}\end{ex}

\end{document}